\documentclass{article}
\usepackage{latexsym,amsfonts,amsmath,amsthm,amssymb,makeidx}
\usepackage[title]{appendix}
\usepackage{CJK,CJKnumb,CJKulem,times,dsfont,ifthen,mathrsfs,latexsym,amsfonts, color}
\usepackage{amsmath,amsthm,makeidx,fontenc,amssymb,bm,graphicx,psfrag,listings, curves,extarrows}

\let\oldbibliography\thebibliography
\renewcommand{\thebibliography}[1]{%
\oldbibliography{#1}%
\setlength{\itemsep}{0pt}%
}

\makeindex
\newtheorem{definition}{Definition}[section]
\newtheorem{theorem}{Theorem}[section]
\newtheorem{lemma}{Lemma}[section]
\newtheorem{corollary}{Corollary}[section]

\newtheorem{remark}{Remark}[section]

\newcommand{\Om}{\Omega}

\newcommand{\dd}{\delta}

\newcommand{\la}{\lambda}
\newcommand{\La}{\Lambda}
\newcommand{\pa}{\partial}

\newcommand{\R}{\mathbb R}

\newcommand{\bb}{\beta}

\newcommand{\bt}{\begin{theorem}}
\newcommand{\et}{\end{theorem}}
\newcommand{\bl}{\begin{lemma}}
\newcommand{\el}{\end{lemma}}
\newcommand{\bd}{\begin{definition}}
\newcommand{\ed}{\end{definition}}
\newcommand{\bc}{\begin{corollary}}
\newcommand{\ec}{\end{corollary}}
\newcommand{\bp}{\begin{proof}}
\newcommand{\ep}{\end{proof}}
\newcommand{\bx}{\begin{example}}
\newcommand{\ex}{\end{example}}
\newcommand{\bi}{\begin{exercise}}
\newcommand{\ei}{\end{exercise}}
\newcommand{\bo}{\begin{prop}}
\newcommand{\eo}{\end{prop}}
\newcommand{\br}{\begin{remark}}
\newcommand{\er}{\end{remark}}
\newcommand{\be}{\begin{equation}}
\newcommand{\ee}{\end{equation}}
\newcommand{\ba}{\begin{align}}
\newcommand{\ea}{\end{align}}
\newcommand{\bn}{\begin{enumerate}}
\newcommand{\en}{\end{enumerate}}
\newcommand{\bg}{\begin{align*}}
\newcommand{\bcs}{\begin{cases}}
\newcommand{\ecs}{\end{cases}}

\newcommand{\pl}{\partial}

\newcommand{\bean}{\begin{eqnarray*}}
\newcommand{\eean}{\end{eqnarray*}}

\numberwithin{equation}{section}

\begin{document}

\title{\bf  Existence, nonexistence, symmetry and uniqueness of  ground state for    critical Schr\"odinger system involving  Hardy term\thanks{Supported by NSFC.
E-mail address:
luosp14@mails.tsinghua.edu.cn(Luo);\,\, wzou@math.tsinghua.edu.cn
(Zou)}}
\date{}
\author{{\bf Senping Luo   \&   Wenming Zou }\\
\footnotesize {\it   Department of Mathematical Sciences, Tsinghua University, Beijing 100084, China}
}

\maketitle
\begin{center}
\begin{minipage}{120mm}
\begin{center}{\bf Abstract}\end{center}
We study the following elliptic system with critical exponent:
\begin{displaymath}
\begin{cases}-\Delta u_j-\frac{\lambda_j}{|x|^2}u_j=u_j^{2^*-1}+\sum\limits_{k\neq j}\beta_{jk}\alpha_{jk}u_j^{\alpha_{jk}-1}u_k^{\alpha_{kj}},\;\;x\in\R^N,\\
u_j\in D^{1,2}(\R^N),\quad  u_j>0 \;\; \hbox{in} \quad \R^N\setminus \{0\},\quad  j=1,...,r.\end{cases}\end{displaymath}
Here $N\geq 3, r\geq2,  2^*=\frac{2N}{N-2},  \lambda_j\in (0, \frac{(N-2)^2}{4})$ for all $ j=1,...,r $; $\beta_{jk}=\beta_{kj}$; \;  $\alpha_{jk}>1, \alpha_{kj}>1,$ satisfying $\alpha_{jk}+\alpha_{kj}=2^* $ for all $k\neq j$.
Note that the nonlinearities $u_j^{2^*-1}$  and the coupling terms all  are  critical in arbitrary dimension $N\geq3 $.
The signs of the coupling constants $\bb_{ij}$'s are decisive for the existence of  the ground state solutions.
 We show that the critical system with $r\geq 3$ has a positive least energy solution for all $\beta_{jk}>0$. However,
  there is no ground state solutions if all $\beta_{jk}$ are negative.
 We also prove  that the positive solutions  of the  system are radially symmetric.
 Furthermore, we obtain  the  uniqueness theorem for the case $r\geq 3$
  with $N=4$ and the existence theorem when  $r=2$ with general coupling exponents.

\vskip0.10in


\end{minipage}
\end{center}

\vskip0.30in
\section{Introduction}

Consider the solitary wave solutions to the time-depending  $r$-coupled nonlinear Schr\"odinger equations given by
 \be\label{2sys1} \begin{cases}-i\frac{\pl}{\pl t}\Phi_j=\Delta\Phi_j-a_{j}(x)\Phi_j+\mu_j|\Phi_j|^2\Phi_j+\sum_{i\neq j}\bb_{ij}|\Phi_i|^2\Phi_j, \\
 \Phi_j=\Phi_j(x,t)\in\mathbb{C},\;\;j=1,2,...,r;  \quad \;x\in \R^N,\; t>0,\\
 \Phi_j(x,t)\rightarrow 0,\;\;\hbox{as}\;\; |x|\rightarrow +\infty, \;t>0, \;j=1,2,...,r.
 \end{cases}
 \ee
Where $\mu_j>0$ are positive constants and $\bb_{ij}$'s are coupling constants;  $a_j(x)$ are potential functions.
When $N\leq 3$,   system \eqref{2sys1} appears in many physical problems, especially in nonlinear optics.
Physically, the solution $\Phi_j$ denotes the $j^{th}$ component of the beam in Kerr-like photorefractive media.
The positive constant  $\mu_j$ is standing  for the  self-focusing in the $j^{th}$ component of the beam.
The coupling constant  $\bb_{ij}$ represents  the interaction between the $i^{th}$ and the $j^{th}$ component of the beam.
As $\bb_{ij}>0$, the interaction is attractive, but the interaction is repulsive if $\bb_{ij}<0$.
To obtain the  solitary wave solutions of the system \eqref{2sys1}, ones usually  set   $\Phi_j(x,t)=e^{i\la_jt}u_j(x)$ and
 may transform the system \eqref{2sys1} to steady-state $r$-coupled nonlinear Schr\"odinger  system:
\be\label{2 1.3}
\begin{cases}-\Delta u_j+V_j(x) u_j=\mu_ju_j^3+\sum\limits_{k\neq j}\beta_{kj}u_k^2u_j,\;\;x\in\R^N,\\
u_j\geq0,\;\; x\in\R^N,\;\;u_j\rightarrow0\;\;\hbox{as}\;\; |x|\rightarrow +\infty;\;\;j=1,2,...,r.\end{cases}\ee
We briefly  recall some previous works on this line.

\vskip0.1in
\noindent{\bf Subcritical case:}
When $N\leq 3$, then the critical Sobolev exponent $2^\ast:=\frac{2N}{N-2}\in [6, +\infty]$ and hence the nonlinear terms (including the coupling terms) of \eqref{2 1.3}  are of
subcritical growth. For such cases, we call the system \eqref{2 1.3}  subcritical which has received great interest in the last decade and large
number of papers published. On this line, although we can not exhaustedly  enumerate and    all those articles,  we refer the readers to  \cite{Ambrosetti2007,Bartsch2006,Bartsch2007, ACR-zou9,Bartsch2010,Byeon2015,chang-lin-lin-lin,Chen-lin-zou=pisa,Chen-zou-optimal,Dancer2010,Eugenio-Benedetta-Squassina,LiuWang-CMP,Maia2006,Noris2010,Lin2005,Lin2005a,Lin2006,peng-wang=ARMA,Pomponio2006,sato-wang=CVPDE,sato-wang=ANS,Sirakov2007,Soave=1,Soave=2,TV-zwm9, Wei2007,Wei2008,Wei-Weth-2008,Wei-Yao-2012} and the references cited therein for various  existence of solutions.

\vskip0.1in
\noindent{\bf Critical case:}  When $N=4$, then the critical Sobolev exponent $2^\ast=4$ and thus the nonlinear terms (including the coupling terms) of \eqref{2 1.3} all are of critical growth. Due to the lack of compactness, this kind problems become thorny. Basically, such a system \eqref{2 1.3} with $r=2$ and $V_j=const$ was firstly studied in  \cite{Chen2012} (including the same system defined on a bounded domain). The positive least energy solutions
and phase separation were obtained in  \cite{Chen2012}.   Later, the higher dimension case (i.e., $N\geq 5$) was also considered in \cite{Chen2015} where some different phenomenon from the  3-D and 4-D cases were  observed.  We also note that, in \cite{wang-willem}, a partial symmetry was involved when $N=4$ (and $N=2,3$) under the premise of assuming the existence of  the minimizer.

\vskip0.13in
In the current paper, we are interested in    the following   $r$-coupling system:

\be\label{model}
\begin{cases}-\Delta u_j-\frac{\lambda_j}{|x|^2}u_j=u_j^{2^*-1}+
\sum\limits_{k\neq j}\beta_{jk}\alpha_{jk}u_j^{\alpha_{jk}-1}u_k^{\alpha_{kj}},\;\;x\in\R^N,\\
u_j\in D^{1,2}(\R^N),\;\;u_j>0\;\;\hbox{in} \;\;\R^N\setminus \{0\},\;\;j=1,...,r;
\end{cases}
\ee
where  $N\geq 3, r\geq2,    \lambda_j\in (0, \frac{(N-2)^2}{4})$ for all $ j=1,...,r $;
 and $\beta_{jk}=\beta_{kj}$, $\alpha_{jk}>1, \alpha_{kj}>1,$ satisfying $\alpha_{jk}+\alpha_{kj}=2^* $ for all $k\neq j$. Note that
  $\alpha_{jk} \not= \alpha_{kj}$ is allowed. We are concerned
with   the existence, nonexistence, symmetry and uniqueness of the  ground state for  the system \eqref{model}.

\vskip0.1in
When $V_j(x)=-\frac{\lambda_j}{|x|^2}$,  the Hardy's type potentials  appear, then   the system \eqref{model}   arises  in
several physical contexts including nonrelativistic quantum mechanics, molecular physics, quantum cosmology, and linearization of combustion models.
 The Hardy's type potentials do not belong to Kato's class, so they cannot be regarded as a lower order perturbation term. In particular, any nontrivial solution   is  singular at $x = 0$.  We refer to the  papers \cite{Abdellaoui2004,Felli2006,Smets2005,Terracini1996} for the scalar
 equations.

\vskip0.12in
For the case of  $r=2$, the two-coupled  system  \eqref{model}   has been studied  in \cite{Chen2015a} where the positive ground state solutions are obtained and are all radially symmetric. It turns out
that the least energy level depends heavily on the relations among  $\alpha_{jk}$ and  $\alpha_{kj}$. Besides, for sufficiently small coupling constants, positive solutions are also obtained via a variational perturbation approach.  It is point out that    the Palais-Smale condition
cannot hold for any positive energy level, which makes the study via variational
methods rather complicated, see \cite{Chen2015a}. We remark that  in \cite{Chen-zou-remark}, when the coupling constant is replaced by a function
decaying  to zero, then the existence of ground state is obtained.

\vskip0.12in

   However, when $r>2$, the study of  system  \eqref{model} becomes rather  complicated.  In particular, even for the two-coupled case of  \eqref{model} (i.e., $r=2$),  the characteristics and uniqueness of the least energy solution to \eqref{model} have not been solved completely in  \cite{Chen2015a} (see Remarks \ref{zwm-rmk=1}-\ref{zwm-rmk=2} below).  In the present paper, we give some positive answers for several standing problems related to the system  \eqref{model}.
   We will introduce some quite different  techniques  than usual. Precisely, we shall study some  nonlinear constraint problems which will
   play an important role for exploring the  multi-coupled system  \eqref{model}. We consummate
   the results due to  \cite{Chen2012, Chen-zou-remark,Chen2015a,Chen2015}.

\vskip0.10in

Let  $\lambda_j \in (0,\Lambda_N)$ for all $j=1,...,r$, where $\Lambda_N:=\frac{(N-2)^2}{4} $.  Set
\be \|u\|_{\la_j}^2:=\int_{\R^N}|\nabla u|^2-\frac{\la_j}{|x|^2}u^2; \quad  \langle u,v\rangle_{\la_j}:=\int_{\R^N}\nabla u\nabla v-\frac{\la_j}{|x|^2}uv,\ee
for all $u, v\in D^{1,2}:=D^{1,2}(\R^N).$ Denote the norm of $L^p(\R^N)$ by $|u|_p=(\int_{\R^N} |u|^p)^{\frac{1}{p}}$. Let
\be\label{zwm-aug-15} I_{\la_j}(u):=\frac{1}{2}\|u\|_{\la_j}^2-\frac{1}{2^*}\int_{\R^N}|u|^{2^*},\quad\quad j=1,...,r.
\ee
 We call a solution $(u_1,...,u_r)$ of \eqref{model} nontrivial if all $u_j\not\equiv 0, j=1,...,r$.
  We call that a solution $(u_1,...,u_r)$ is positive if all $u_j>0$ in $\R^N\setminus \{0\}$ for all $j=1,...,r$.
   We call a solution $(u_1,...,u_r)\not=(0,..., 0)$ is semi-trivial if there exists  some $i_0$ satisfying  $u_{i_0}\equiv 0$.
Throughout this paper, we are only interested in nontrivial solutions of \eqref{model}. Define $\mathbb{D}:=D^{1,2}\times \cdots \times D^{1,2}$ with the norm
\be\|(u_1,...,u_r)\|_{\mathbb{D}}^2:=\sum_{j=1}^r\|u_j\|_{\la_j}^2.\ee
Then the nontrivial solutions of \eqref{model}  correspond  to the    nontrivial critical points of the $C^1$ functional $J:\mathbb{D}\rightarrow \R,$ where
\be\label{2JJJ} J(u_1,...,u_r):=\sum_{j=1}^r I_{\la_j}(u_j)-\frac{1}{2}\sum_{1\leq j\neq k \leq r} \beta_{jk}\int_{\R^N}|u_j|^{\alpha_{jk}}|u_k|^{\alpha_{kj}}.
\ee

\bd We say a solution $(u_{1,0},...,u_{r,0})$ of \eqref{model} is a ground state solution if $(u_{1,0},...,u_{r,0})$ is nontrivial and $J(u_{1,0},...,u_{r,0})\leq J(u_1,...,u_r)$ for any other nontrivial solution $(u_1,...,u_r)$ of \eqref{model}.
\ed

To obtain the ground state solutions of \eqref{model}, we define the Nehari manifold:
\be\label{zwm=110}
 \aligned
 \mathcal{N}:=\{&(u_1,...,u_r)\in \mathbb{D}: \quad   u_j\not \equiv 0  \hbox{ and }\\
 &\quad \|u_j\|^2_{\la_j}=\int_{\R^N}|u_j|^{2^*}+\sum_{k\neq j} \beta_{jk}\alpha_{jk}|u_j|^{\alpha_{jk}}|u_k|^{\alpha_{kj}},\;\; j=1,...,r\}.
 \endaligned
\ee
Then any nontrivial solution of \eqref{model}  belongs  to $\mathcal{N}$. Note that $\mathcal{N}\neq \emptyset.$ We set
\be\label{zwm=111} \Theta:=\inf_{(u_1,...,u_r)\in \mathcal{N}}J(u_1,...,u_r).  \ee
Hence, $\Theta=\inf_{(u_1,...,u_r)\in \mathcal{N}} \frac{1}{N}\sum_{j=1}^r\|u_j\|_{\la_j}^2.$
It is easy to see  that $ \Theta>0$.

\vskip0.1in
 Recall the following scalar  equation which
has been deeply investigated in the literature (see for example \cite{Terracini1996}):
\be\label{2hardyeq}\begin{cases}-\Delta u-\frac{\lambda_i}{|x|^2}u=u^{2^*-1},\;\; x\in\R^N,\\
u\in D^{1,2}(\R^N),\;\;u> 0\;\;\hbox{in}\;\;\R^N\setminus \{0\},
\end{cases}\ee
which has exactly an  one-dimensional $C^2$-manifold of positive solutions given by
 \be\label{2zhardy} Z_i:=\{z_{\mu}^i(x)=\mu^{-\frac{N-2}{2}}z_1^i(\frac{x}{\mu}):    \;\;\mu> 0\},   \ee
 where
 \be z_1^i(x)=\frac{A(N,\la_i)}{{|x|^{a_{\la_i}}(1+|x|^{2-\frac{4a_{\la_i}}{N-2}})^{\frac{N-2}{2}}}}  \nonumber \ee
 and
 $$a_{\la_i}=\frac{N-2}{2}-\sqrt{\frac{(N-2)^2}{4}-\la_i},\quad A(N,\la_i)=\frac{N(N-2-2a_{\la_i})^2}{N-2}.$$
  Moreover, all positive solutions of $\eqref{2hardyeq}$
 satisfy
\be\label{zwm=116}I_{\la_i}(z_{\mu}^i)=\frac{1}{N}\|z_{\mu}^i\|_{\la_i}^2=\frac{1}{N}|z_{\mu}^i|^2_{2^*}=\frac{1}{N}S(\la_i)^{\frac{N}{2}},
\ee
where
\be \label{2ha_in}S(\la_i):=\inf_{u\in D^{1,2}(\R^N)\setminus \{0\}} \frac{\|u\|_{\la_i}^2}{|u|^2_{2^*}}
 =\frac{\|z_{\mu}^i\|_{\la_i}^2}{|z_{\mu}^i|^2_{2^*}}=\Big(1-\frac{4\la_i}{(N-2)^2}\Big)^{\frac{N-1}{N}}S
\nonumber \ee
and  $S$ is the sharp constant of $D^{1,2}(\R^N)\hookrightarrow L^{2^*}(\R^N)$ (see e.g.,  \cite{Talenti1976}):
\be \int_{\R^N}|\nabla u|^2 \geq S\Big(\int_{\R^N}|u|^{2^*}\Big)^{\frac{2}{2^*}}. \nonumber
\ee
In the current paper, we always assume that $\beta_{jk}=\beta_{kj}$ for $1\leq j,k\leq r$.
Now we are ready to state the main theorems of this article.

\bt\label{2th1}Assume that $N\geq 3, \la_j \in (0,\Lambda_N)$ and $\alpha_{jk}>1,\alpha_{kj}>1, \alpha_{jk}+\alpha_{kj}=2^* $ for  $1\leq j,k\leq r$.
\begin{itemize}
\item[(1)] (Nonexistence) If $\beta_{jk}<0, \forall k\neq  j$, then $\displaystyle  \Theta\equiv \sum_{j=1}^N\frac{1}{N}\Big((1-\frac{4\la_j}{(N-2)^2})^{\frac{N-1}{N}}S\Big)
^{\frac{N}{2}}$, and $\Theta$ cannot be attained, i.e.,  there is no ground state solution to \eqref{model}.
\item[(2)] (Existence) Let $\beta_{jk}>0,\forall k\neq  j,$ satisfy
$$\Big(r+\sum_{j,k=1,j\neq k}^r\frac{2^*}{2}\beta_{jk}\Big)\Big/\Big(\max_{j,l}B_{j,l}\Big)>\Big(1+\sum_{\alpha=1}^{r-1}\frac{\La_N-\la_{\alpha}}{\La_N-\la_r}
\Big)^{\frac{N}{N-2}},$$
then \eqref{model} has a positive ground state solution $(u_1,...,u_r)\in \mathbb{D}$, which is radially symmetric and whose
energy  satisfies $$\Theta<\min_{l\in\{1,2,...,r\}}\min_{j\in A^l}B_{j,l}^{-\frac{N-2}{2}}\frac{1}{N}S(\la_j)^{\frac{N}{2}},$$
where $B_{j,l}=\sum_{k\neq j,k\in A^l}\beta_{jk}\alpha_{jk}+1$,\;\; $A^l=\{1,2,...,r\}\setminus\{l\}$.
\end{itemize}
\et


\vskip0.36in

\bt\label{2th2} Assume that $N=3$ or $N=4$, $\alpha_{jk}+\alpha_{kj}=2^*,\alpha_{jk}\geq2,\alpha_{kj}\geq2,\la_j\in(0,\Lambda_N),\beta_{jk}>0$,
  $\forall k\neq j$, then any positive solution of \eqref{model} is radially symmetric with respect to the origin.
\et

\vskip0.3in
Next, we   obtain the existence and uniqueness results about  the ground state
 to  the following critical elliptic system in $\R^4$ involving the  Hardy's singular  term:
\be\label{2m2}\begin{cases}-\Delta u_j-\frac{\lambda}{|x|^2} u_j
=\gamma_{jj}u_j^3+\sum\limits_{i\neq j}\gamma_{ij}u_i^2u_j,\\
\quad u_j(x)>0, \quad j=1,...,r, \quad \quad x\in \R^4\setminus\{0\}.\end{cases}\ee
We have the following result.
\bt\label{2th3} Considering the system \eqref{2m2}.  Assume that $$N=4,\;  r\geq3,\;\la\in(0, \Lambda_N), \;  \gamma_{ji}=\gamma_{ij}, \;\; \det(\gamma_{ij})\neq 0, \;\;  \sum_{k}\gamma^{kj}>0,\;\; i, j=1,...,r;$$ where the matrix $(\gamma^{jk})$ represents the inverse matrix of $(\gamma_{ml})$. Then
\begin{itemize}
\item[(1)](existence) $(\sqrt{c_1}z_{\mu}^1,...,\sqrt{c_r}z_{\mu}^1) $ ($ \mu>0$) is a positive least energy solution of \eqref{2m2}, where $z_{\mu}^1$ is a solution of (see \eqref{2zhardy})
\be\label{zwm-aug-th13}\begin{cases}-\Delta u-\frac{\lambda}{|x|^2}u=u^{3},\;\;x\in\R^N,\\
u\in D^{1,2}(\R^N),\quad u> 0\;\;\hbox{in} \;\;\R^N\setminus \{0\},
\end{cases}\ee
and the constant ${c_j}>0$ satisfying
\be  \sum_{k=1}\gamma_{jk}c_k=1, \quad \;\; j=1,...,r. \nonumber \ee

\item[(2)](uniqueness) let $(u_1,u_2,...,u_r)$ be any   least energy solution of \eqref{2m2}, then $(u_1,u_2,...,u_r)$$=(\sqrt{c_1}z_{\mu}^1$
$,...,\sqrt{c_j}z_{\mu}^1)$, where
\end{itemize}
\be  \sum_{k=1}\gamma_{jk}c_k=1, \;\; j=1,...,r. \nonumber \ee
\et

\vskip0.10in
\br\label{zwm-rmk=1} When $r=2, \la=0$, the existence of the ground state for   system \eqref{2m2} in $\R^4$ was firstly studied in \cite{Chen2012}.
\er


\vskip0.3in

 Lastly, we consider the following two-coupled  doubly critical shr\"{o}dinger system:
 \begin{equation}\label{2mod4}
 \begin{cases}
 -\Delta u-\frac{\lambda}{|x|^2}u=u^{2^*-1}+\nu \alpha u^{\alpha-1}v^{\beta},  \;\; x\in \R^N,  \\
 -\Delta v-\frac{\lambda}{|x|^2}v=v^{2^*-1}+\nu \beta  u^{\alpha}v^{\beta-1},  \;\; x\in \R^N.
 \end{cases}
 \end{equation}
We have the following theorem.

\bt\label{2th4} In the system \eqref{2mod4}, we assume that  $\la\in(0,\Lambda_N),1<\alpha,\beta<2$  and $\alpha+\beta=2^*$ (these imply
 $N\geq 5$).
 \begin{itemize}
\item [(1)]  If $\nu>0$, then   $(c_1 z_{\mu}, c_2 z_{\mu})$  is a positive solution of \eqref{2mod4} for any $\mu>0$, where $z_{\mu}$ is a solution of  the following equation:
    \be\label{zou=999}\begin{cases}-\Delta u-\frac{\lambda}{|x|^2}u=u^{2^*-1},\;\; x\in\R^N,\\
u\in D^{1,2}(\R^N),\;\;u> 0\;\;\hbox{in}\;\;\R^N\setminus \{0\},
\end{cases}\ee

\item [(2)]  If
\be\nu>(\frac{2^\ast}{2}-1)/\min\{d_1(\alpha,\beta),d_2(\alpha,\beta), d_3(\alpha,\beta)\}, \nonumber\ee then
$(c_1 z_{\mu},c_2 z_{\mu})$ is a positive ground state solution of \eqref{2mod4}, where $c_1,c_2$ are the roots of the algebraic system about $(x_1, x_2)$:
\begin{equation}\nonumber
\begin{cases}
&x_1^{\frac{2^\ast}{2}-1}+\nu\alpha x_1^{\frac{\alpha}{2}-1}x_2^{\frac{\beta}{2}}=1, \\
&x_2^{\frac{2^\ast}{2}-1}+\nu\beta x_1^{\frac{\alpha}{2}}x_2^{\frac{\beta}{2}-1}=1,
\end{cases}
\end{equation}
and $d_1(\alpha,\beta),d_2(\alpha,\beta),d_3(\alpha,\beta)$ are defined as following:
\be\nonumber d_1(\alpha,\beta)=2^\ast(1-\frac{\alpha}{2})^{\frac{\alpha}{2^\ast}}(1-\frac{\beta}{2})^{\frac{\beta}{2^\ast}}\;\;\hbox{if}\;\;\alpha\neq\beta;
\quad d_1(\alpha,\beta)=2^\ast/2\;\;\hbox{if}\;\;\alpha=\beta; \nonumber
\ee
\be\nonumber d_2(\alpha,\beta)=\beta(1-\frac{\beta}{2})^{1-\frac{\alpha}{2}}(1-\frac{\alpha}{2})^{\frac{\alpha}{2}}
+\frac{1}{2}\alpha\beta(1-\frac{\alpha}{2})^{1-\frac{\alpha}{2}}(1-\frac{\beta}{2})^{\frac{\alpha}{2}-1};\quad\quad\quad\nonumber
\ee
\be\nonumber d_3(\alpha,\beta)=\alpha(1-\frac{\alpha}{2})^{1-\frac{\beta}{2}}(1-\frac{\beta}{2})^{\frac{\beta}{2}}
+\frac{1}{2}\alpha\beta(1-\frac{\alpha}{2})^{1-\frac{\beta}{2}}(1-\frac{\beta}{2})^{\frac{\beta}{2}-1}.\quad\quad\quad \nonumber
\ee
\end{itemize}

\et


\br\label{zwm-rmk=2} When $N\geq 5$, the existence of ground state solution is essentially proved in \cite{Chen2015a}. Here, the further characteristics
 is given. If $\alpha=\beta$, then
$$(\frac{2^\ast}{2}-1)/\min\{d_1(\alpha,\beta),d_2(\alpha,\beta),d_3(\alpha,\beta)\}=\frac{2}{N}.$$
We remark that, for the special case of  $\alpha=\beta=\frac{2^\ast}{2}$, the uniqueness of the
ground state solution of \eqref{2mod4} was   obtained by a different method in \cite{Chen2015a}.
\er

 The paper is organized as follows.  In Section 2,
 we develop several lemmas which will also have other applications.
We give the the proof of Theorem \ref{2th1} in Section 3, where we will use the concentration-compactness
 principle due to  \cite{Lions1985,Lions1985a}. In Section 4, Theorem \ref{2th2} is proved by the moving plane method.
  In Section 5, we firstly construct some powerful lemmas  and then
 obtain  the existence and uniqueness results about the positive ground state.  Theorems \ref{2th3}-\ref{2th4} will get  proved there.

\section{Preliminaries}


 We firstly  deal with  the following  nonlinear  algebraic equations
 which is important for construct the nonexistence of the  ground state solution.


\bl\label{2lm2}  Assume
\be\label{zwm-aug-1} C_j:=|u_j|_{2^*}^{2^*}-\sum_{k\neq j}\frac{\alpha_{jk}^2-\alpha_{jk}\alpha_{kj}}{2^*}\int_{\R^N}|\beta_{jk}||u_j|^{\alpha_{jk}}|u_k|^{\alpha_{kj}}>0,\;\; j=1,...,r.
\ee
Consider  the algebraic equations about $t_j$:
\be\label{2ageq} t_j^2\|u_j\|^2_{\la_j}=t_j^{2^*}|u_j|_{2^*}^{2^*}+\sum_{k\neq j}t_j^{\alpha_{jk}}t_k^{\alpha_{kj}}\int_{\R^N}\beta_{jk}\alpha_{jk}|u_j|^{\alpha_{jk}}|u_k|^{\alpha_{kj}},\;\; j=1,...,r,\ee
where $\beta_{jk}<0,\alpha_{jk}>1,\alpha_{kj}>1,\alpha_{jk}+\alpha_{kj}=2^*,  u_j\not\equiv 0.$
Then we have the following priori estimate for the positive solution of \eqref{2ageq} (if any):
\be\label{zwm-aug-11}\min_{1\leq j\leq r}\Big(\frac{A_j}{B_j}\Big)^{\frac{1}{\alpha}} \leq \Big(\frac{\|u_j\|_{\la_j}^2}{|u_j|^{2^*}_{2^*}}\Big)^{\frac{1}{\alpha}}\leq t_j \leq \Big(\frac{\sum_{j=1}^rA_j}{\min\{C_1,...,C_r\} } \Big)^{\frac{1}{\alpha}},
\ee
where
 \be\alpha=2^*-2=\frac{4}{N-2},\; \; A_j=\|u_j\|^2_{\la_j}>0,\;\;  B_j=|u_j|_{2^*}^{2^*}>0,\;\; j=1,...,r. \nonumber\ee
In particular, the systems  \eqref{2ageq} has  a positive solution   provided that
\be\label{2con2} d:=\frac{1}{\alpha}r\max_j\Big(\frac{A_j}{B_j}\Big)^{\frac{1}{\alpha}}\Big(1+\max_j \max f_j\Big)\max_{j,m}\max|\frac{\pa f_j}{\pa t_m}|<1,
\ee
where
\be\label{2fj} f_j(t_1,...,t_r):=\frac{1}{A_j}\sum_{k\neq j}t_j^{\alpha_{jk}-2}t_k^{\alpha_{kj}}\int_{\R^N}|\beta_{jk}|\alpha_{jk}|u_j|^{\alpha_{jk}}|u_k|^{\alpha_{kj}}.
\ee
Furthermore, if
\be\label{2bbar}\bar{\beta}_{jk}:=\int_{\R^N}|\beta_{jk}||u_j|^{\alpha_{jk}}|u_k|^{\alpha_{kj}},\;\; \forall j \neq k,\ee
all are small enough, then \eqref{2ageq} admits a positive solution $(t_1,...,t_r)$. In particular, each $t_s$ ($s=1, ..., r$)  satisfying
\be t_s\rightarrow \Big(\frac{\|u_s\|_{\la_s}^2}{|u_s|^{2^*}_{2^*}}\Big)^{\frac{1}{\alpha}}\;\;\hbox{as all  } \bar{\beta}_{jk}\to 0, \;\forall j\not=k.
\ee
\el

\br In Lemma \ref{2lm2} above, $\max f_j=\max f_j(t_1,...,t_r)$ which is a finite value in view of the priori bound
on $(t_1,...,t_r)$ obtained in \eqref{zwm-aug-11}. The same conclusion is true for $\max|\frac{\pa f_j}{\pa t_m}|$.

\er


\bp
Firstly, by definition we have $\min_{1\leq j\leq r}\Big(\frac{A_j}{B_j}\Big)^{\frac{1}{\alpha}} \leq \Big(\frac{\|u_j\|_{\la_j}^2}{|u_j|^{2^*}_{2^*}}\Big)^{\frac{1}{\alpha}}$.
Since $\beta_{jk}<0, \forall j \neq k$, then  we have
\be t_j^2\|u_j\|^2_{\la_j}\leq t_j^{2^*}|u_j|_{2^*}^{2^*}, \nonumber
\ee
that is,
\be t_j\geq(\frac{\|u_j\|_{\la_j}^2}{|u_j|^{2^*}_{2^*}})^{\frac{1}{\alpha}}. \nonumber
\ee
Recall that
$ \alpha_{jk}+\alpha_{kj}=2^*$, by Young's inequality, we have
$$ t_j^{\alpha_{jk}}t_k^{\alpha_{kj}}\leq \frac{\alpha_{jk}}{2^*}t_j^{2^*}+\frac{\alpha_{kj}}{2^*}t_k^{2^*},
$$
then
\begin{equation}\label{2inq1}
 \aligned
t_j^2A_j&=t_j^{2^*}B_j-\sum_{k\neq j}\int_{\R^N}|\beta_{jk}|\alpha_{jk}|u_j|^{\alpha_{jk}}|u_k|^{\alpha_{kj}}t_j^{\alpha_{jk}}t_k^{\alpha_{kj}}\\
&\geq B_j t_j^{2^*}-\sum_{k\neq j}D_{jk}\frac{\alpha_{jk}}{2^*}t_j^{2^*}+D_{jk}\frac{\alpha_{kj}}{2^*}t_k^{2^*}\\
&=(B_j-\sum_{k\neq j}D_{jk}\frac{\alpha_{jk}}{2^*})t_j^{2^*}-\sum_{k\neq j} D_{jk}\frac{\alpha_{kj}}{2^*}t_k^{2^*},
 \endaligned\end{equation}
where
\be D_{jk}:=\int_{\R^N}|\beta_{jk}|\alpha_{jk}|u_j|^{\alpha_{jk}}|u_k|^{\alpha_{kj}}. \nonumber
\ee
Summing up  \eqref{2inq1} from $j=1$ to $j=r$, thus
\begin{equation}
 \aligned
 \sum_{j=1}^r A_j t_j^2&\geq \sum_{j=1}^r(B_j-\sum_{k\neq j} D_{jk}\frac{\alpha_{jk}}{2^*})t_j^{2^*}-\sum_{j=1}^r\sum_{k\neq j}D_{jk}\frac{\alpha_{kj}}{2^*})t_k^{2^*}\\
 &=\sum_{j=1}^r(B_j-\sum_{k\neq j} D_{jk}\frac{\alpha_{jk}}{2^*})t_j^{2^*}-\sum_{k=1}^r\sum_{j\neq k}D_{jk}\frac{\alpha_{kj}}{2^*})t_k^{2^*}\\
 &=\sum_{j=1}^r(B_j-\sum_{k\neq j} D_{jk}\frac{\alpha_{jk}}{2^*})t_j^{2^*}-\sum_{j=1}^r\sum_{k\neq j}D_{kj}\frac{\alpha_{jk}}{2^*})t_j^{2^*}\\
 &=\sum_{j=1}^r(B_j-\sum_{k\neq j}\frac{\alpha_{jk}(D_{jk}-D_{kj})}{2^*})t_j^{2^*}\\
 &=\sum_{j=1}^r C_j t_j^{2^*}.
 \endaligned\end{equation}
 Recall that $A_j>0,C_j>0,j=1,...,r$.  For the positive solution of \eqref{2ageq}, without loss of generality, we assume that $t_1=\max\{t_1,...,t_t\}$, then we have
 \be C_1 t_1^{2^*}\leq \sum_{j=1}^r A_j t_j^2\leq \sum_{j=1}^r A_j t_1^2, \nonumber
 \ee
that is,
\be t_1:=\max_j t_j\leq \Big(\frac{1}{C_1}\sum_{j=1}^r A_j\Big)^{\frac{1}{\alpha}}\leq \Big(\frac{1}{\min_{j}{C_j}}\sum_{j=1}^r A_j\Big)^{\frac{1}{\alpha}}.
\ee
Hence the priori estimate is obtained. Hence, there are two positive constants $T_1>0, T_2>0$ such that for all   positive solution $t_j$  of \eqref{2ageq}:
$$  t_j\in [T_1, T_2], \;\;\; \forall j=1,..., r.$$
In the following, we will use Picard's iteration to obtain the existence of positive solution of \eqref{2ageq}.
Recall the notation of $f_j$ in \eqref{2fj}, the equation \eqref{2ageq} becomes
\begin{equation}
 \aligned
t_j&=\Big(\frac{A_j}{B_j}+\frac{A_j}{B_j}f_j(t_1,...,t_r)\Big)^{\frac{1}{\alpha}}\\
&=\Big(\frac{A_j}{B_j}\Big)^{\frac{1}{\alpha}}\Big(1+f_j(t_1,...,t_r)\Big)^{\frac{1}{\alpha}}, \;\; j=1,...,r.
 \endaligned\end{equation}
We   select arbitrarily  an initial value   $t_0=(t_{1,0},...,t_{r,0})\in [T_1, T_2]^r$, then
\begin{equation}
 \aligned
&|t_{j,n+1}-t_{j,n}|\\
&=\Big(\frac{A_j}{B_j}\Big)^{\frac{1}{\alpha}}\Big[\Big(1+f_j(t_{1,n},...,t_{r,n})\Big)^{\frac{1}{\alpha}}
-\Big(1+f_j(t_{1,n-1},...,t_{r,n-1})\Big)^{\frac{1}{\alpha}}\Big]\\
&=\Big(\frac{A_j}{B_j}\Big)^{\frac{1}{\alpha}}\frac{1}{\alpha}(1+f_j(\xi))f_j(\xi)^{\frac{1}{\alpha}-1}\sum_{m=1}^r\frac{\pa f_j}{\pa t_m}(\xi)|t_{m,n}-t_{m,n-1}|\\
&\leq \frac{1}{\alpha}\Big(\frac{A_j}{B_j}\Big)^{\frac{1}{\alpha}}(1+\max_j \max f_j)\max_{j,m}\max |\frac{\pa f_j}{\pa t_m}|\sum_m |t_{m,n}-t_{m,n-1}|,
 \endaligned\end{equation}
where $\xi$ is a vector between $t_n=(t_{1,n},...,t_{r,n})$ and $t_{n-1}=(t_{1,n-1},...,t_{r,n-1})$. Add up  the above inequalities  from $j=1$ to $j=r$, we get that
\begin{equation}
 \aligned
& \sum_{j=1}^r|t_{j,n+1}-t_{j,n}|\\
&\leq r\frac{1}{\alpha}\max_{j}\Big(\frac{A_j}{B_j}\Big)^{\frac{1}{\alpha}}(1+\max_j \max f_j)\max_{j,m}\max |\frac{\pa f_j}{\pa t_m}|\sum_m |t_{m,n}-t_{m,n-1}|\\
&:=d \sum_m |t_{m,n}-t_{m,n-1}|.
 \endaligned\end{equation}
By the assumption \eqref{2con2}, $0<d<1$,  thus we may apply  the classical contraction mapping principle and  know  that the vector sequence  $t_n=(t_{1,n},...,t_{r,n})$ is convergent,  say $t_n=(t_{1,n},...,t_{r,n})\rightarrow t=(t_1,...,t_r)$ as $n\rightarrow \infty$ and  $t$ is  a solution of \eqref{2ageq}. Further,  by our priori estimate,
\be
\min_{j}(\frac{A_j}{B_j})^{\frac{1}{\alpha}} \leq t_{j,n} \leq \Big(\frac{1}{\min\{C_1,...,C_r\}\sum_{j=1}^rA_j}\Big)^{\frac{1}{\alpha}},\;\;  j=1,...,r.
\nonumber \ee
 Let $n\rightarrow \infty $, we have
\be\min_{j}\Big(\frac{A_j}{B_j}\Big)^{\frac{1}{\alpha}} \leq t_j \leq \Big(\frac{1}{\min\{C_1,...,C_r\}\sum_{j=1}^rA_j}\Big)^{\frac{1}{\alpha}}, \;\; j=1,...,r,
\nonumber \ee
it implies that  $t=(t_1,...,t_r)$  is  a  positive solution of \eqref{2ageq}.
Furthermore,  if $\bar{\beta}_{jk}$ (which is defined in \eqref{2bbar}) all are small enough, then it is to see that the solvability conditions \eqref{zwm-aug-1}  and \eqref{2con2} hold. Hence, there exists a positive solution of \eqref{2ageq}. By the priori  estimate of this positive solution and in view of \eqref{2ageq}, we get that $t_j\rightarrow \Big(\frac{\|u_j\|_{\la_j}^2}{|u_j|^{2^*}_{2^*}}\Big)^{\frac{1}{\alpha}},  j=1,...,r$.
\ep

\vskip0.3in

\bl\label{2lm1} If $\Theta$ (which is defined in \eqref{zwm=111}) is attained by  $(u_1,...,u_r)\in \mathcal{N}$, then it is a critical point of $J$ (which is introduced  in \eqref{2JJJ})  provided that  $\beta_{jk}<0,\forall k\neq j$.
\el
\bp Suppose  $\beta_{jk}<0, \forall k\neq j$. Assume that $(u_1,...,u_r)\in \mathcal{N}$ such that $\Theta=J(u_1,...,u_r)$. Define
\be G_j(u_1,...,u_r)=\|u_j\|^2_{\la_j}-\int_{\R^N}\Big(|u_j|^{2^*}+\sum_{k\neq j}\beta_{jk}\alpha_{jk}|u_j|^{\alpha_{jk}}|u_k|^{\alpha_{kj}}\Big),j=1,...,r. \nonumber\ee
Then there exist $L_j\in\R\,  (j=1,...,r)$ such that
\be\label{2Jp} J'(u_1,...,u_r)+\sum_{j=1}^N L_j G_j'(u_1,...,u_r)=0.\ee
Testing \eqref{2Jp} with $(0,...,0,u_i,0,...,0)(i=1,...,r)$, we conclude from $(u_1,...,u_r)\in \mathcal{N}$ that
 \begin{equation}
 \aligned
\big \langle &G'_j,(0,...,u_j,...,0)\big\rangle =2\int_{\R^N}|\nabla u_j|^2-2\int_{\R^N}\frac{\la_j}{|x|^2}u_j^2 \nonumber\\
 &   \quad  -2^*\int_{\R^N}|u_j|^{2^*}-\sum_{k\neq j}\int_{\R^N}\beta_{jk}\alpha_{jk}^2|u_j|^{\alpha_{jk}}|u_k|^{\alpha_{kj}}\nonumber\\
 &=(2-2^*)\int_{\R^N}|u_j|^{2^*}-2\sum_{k\neq j}\int_{\R^N}\beta_{jk}\alpha_{jk}(2-\alpha_{jk})|u_j|^{\alpha_{jk}}|u_k|^{\alpha_{kj}}.\nonumber
 \endaligned\end{equation}
 For $k\neq j$, we have
 \begin{equation}
 \aligned
 \langle &G'_j,(0,...,u_k,...,0)\rangle =-\beta_{jk}\alpha_{jk}\alpha_{kj}\int_{\R^N}|u_j|^{\alpha_{jk}}|u_k|^{\alpha_{kj}}.\nonumber
 \endaligned\end{equation}
Then
\be \sum_{j=1}^N L_j\langle G'_j,(0,...,u_i,...,0)\rangle=0. \nonumber
\ee
Hence
\begin{equation}\label{2lineareq}
 \aligned
\Big((2^*-2)\int_{\R^N}|u_j|^{2^*}&+\sum_{k\neq j}|\beta_{jk}|\alpha_{jk}(2-\alpha_{jk})|u_j|^{\alpha_{jk}}|u_k|^{\alpha_{kj}}\Big)L_j\\
&-\sum_{i\neq j} L_i|\beta_{ji}|\alpha_{ji}\alpha_{ij}\int_{\R^N}|u_j|^{\alpha_{jk}}|u_k|^{\alpha_{kj}}=0.
 \endaligned\end{equation}
Since $(u_1,...,u_r)\in \mathcal{N}$, we have
\be \int_{\R^N}|u_j|^{2^*}>\sum_{k\neq j}|\beta_{jk}|\alpha_{jk}\int_{\R^N}|u_j|^{\alpha_{jk}}|u_k|^{\alpha_{kj}}, \nonumber\ee
hence
\begin{equation}
 \aligned
\Big((2^*-2)&\int_{\R^N}|u_j|^{2^*}+\sum_{k\neq j}|\beta_{jk}|\alpha_{jk}(2-\alpha_{jk})|u_j|^{\alpha_{jk}}|u_k|^{\alpha_{kj}}\Big)\\
&>\sum_{k\neq j}|\beta_{jk}|\alpha_{jk}(2-\alpha_{jk})\int_{\R^N}|u_j|^{\alpha_{jk}}|u_k|^{\alpha_{kj}} \\
&=\sum_{k\neq j}|\beta_{jk}|\alpha_{jk}\alpha_{kj}\int_{\R^N}|u_j|^{\alpha_{jk}}|u_k|^{\alpha_{kj}},\quad j=1,...,r.
 \endaligned\end{equation}
These inequalities above illustrate that the coefficient matrix of \eqref{2lineareq} is diagonally dominant, hence the determinant greater than 0.
Combine with \eqref{2lineareq}, we deduce that $L_j=0 \;(j=1,...,r)$ and then  $J'(u_1,...,u_r)=0$.
\ep


The next  two lemmas   are indispensable  for the  construction of the ground state solution and for  the proof of its uniqueness.
We also believe that they can be applied to  other problems.
\bl\label{2principle} Consider the following  $r+1$ inequalities,
\be \nonumber \begin{cases}
\sum_{k=1}^r x_k\leq 0,\\
f_j(x_1,...,x_r)\geq 0, j=1,...,r,
\end{cases}\ee\\
where $f_j(x_1,...,x_r)$ are nonnegative differentiable functions    with $f_j(0,...,0)=0, j=1,...,r.$    Assume that  the following conditions hold:
\be \det\Big(\frac{\pa f_j}{\pa x_i}\Big)\neq 0;\quad  \sum_{s=1}^r g^{si}>0, \quad  i=1,...,r,\ee
 where $(g_{ij}):=(\frac{\pa f_j}{\pa x_i})$;   the matrix $(g^{kl})$ represents  the inverse matrix of $(g_{ij})$. Then we must have
$ x_j=0$ for all $ j=1,...,r.$

\el

\bp Denote $f_j(x_1,...,x_r)=y_j$, then $y_j\geq 0$ for all $ j=1,...,r$.  Note that
\be\nonumber\sum_{k}{\frac{\pa f_j}{\pa x_k}} \frac{\pa x_k}{\pa y_i}=\delta_{ij},
\ee
that is,
\be\nonumber\sum_{k}g_{jk}\frac{\pa x_k}{\pa y_i}=\delta_{ij}.
\ee
Multiply $g^{sj}$ and sum up  for $j$ in the above equations, we get that
\be\nonumber\sum_j\sum_k g^{sj}g_{jk}\frac{\pa x_k}{\pa y_i}=\sum_j\delta_{ij}g^{sj},
\ee
hence
\be\nonumber\sum_k\delta_{sk}\frac{\pa x_k}{\pa y_i}=g^{si},
\ee
thus we obtain that
\be\nonumber \frac{\pa x_s}{\pa y_i}=g^{si}.
\ee
Then
\be\nonumber \frac{\pa}{\pa y_i}\Big(\sum_{s=1}^r x_s(y_1,...,y_r)\Big)=\sum_{s=1}^r g^{si}>0.
\ee
This means that the function $\sum_{s=1}^r x_s(y_1,...,y_r)$ is strictly increasing in any direction. On the other hand ,
since $y_j\geq 0 $ and $f_j(0,...,0)=0$ for all $j=1,...,r$, combining  with $\sum_{k=1}^r x_k\geq0$,
we obtain
\be\nonumber 0\geq\sum_{s=1}^r x_s(y_1,...,y_r)\geq \sum_{s=1}^r x_s(0,...,0)=0,\ee
it follows that $y_j=0$ and hence  $ x_j=0$ for all $j=1,...r.$
\ep
If $f_j$ does  not satisfy the initial condition $f_j(0,...,0)=0$, then we have the following more  general version than  Lemma \ref{2principle}.


\vskip0.23in

\bl\label{2pringe}
Consider  the following nonlinear constraint problem
\be\nonumber\begin{cases}
\sum_{k=1}^r x_k\leq \sum_{k=1}^r c_k,\\
f_j(x_1,...,x_r)\geq f_j(c_1,...,c_r),j=1,...,r,
\end{cases}\ee
where $f_j(x_1,...,x_r)$ are nonnegative differentiable functions.
Assume the following conditions hold:
\be\label{2zzz2} \det(\frac{\pa f_j}{\pa x_i})\neq 0;\;\;  \sum_{s=1}^r g^{si}>0,\;\;  i=1,...,r,\ee
where $(g_{ij})=(\frac{\pa f_j}{\pa x_i})$,  $(g^{kl})$ represents the inverse matrix of $(g_{ij})$.
Then we must have
$ x_j=c_j$ for all $ j=1,...,r.$
\el

\bp
Take
$ h_j(x_1,...,x_r)=f_j(x_1,...,x_r)-f_j(c_1,...,c_r)$
 and make the transformation $y_j=x_j-c_j,j=1,...,r$.  Let
$ l_j(y_1,...,y_r):=h_j(x_1,...,x_r)$. We may apply  Lemma \ref{2principle}
to  $l_j(y_1,...,y_r) $, then the conclusion follows.
\ep


\vskip0.23in

\br  The  condition   \eqref{2zzz2} in Lemma \ref{2pringe} may be replaced by
\be\nonumber \det(\frac{\pa f_j}{\pa x_i})\neq 0;\quad  \sum_{s=1}^r g^{si}\geq0 , i=1,...,r; \quad \sum_{s=1}^r g^{si_0}>0 \;\hbox{for some}\;\; i_0.\ee
Then the same conclusion as  that  in Lemma \ref{2pringe} holds.
\er


\vskip0.33in

\section{Proof of Theorem \ref{2th1} }

\noindent {\bf The proof of Theorem \ref{2th1}-(1)}. Note  the assumption $\beta_{jk}<0, j\neq k $.  Recall  \eqref{2zhardy}, it is easy to see that $z_{\mu^j}^j\rightharpoonup 0$  weakly in $D^{1,2}(\R^N)$ and so $(z_{\mu^j}^j)^\beta\rightharpoonup 0$  weakly in $L^{2^*/\beta}(\R^N)$ as $\mu\rightarrow\infty$, (here we regard $\mu^j$ as an integer in the expression $z_{\mu}^j$) hence
\be\aligned  &\lim_{\mu\rightarrow +\infty}|\beta_{jk}| \int_{\R^N}(z_{\mu^i}^i)^{\alpha}(z_{\mu^j}^j)^{\beta}dx\\
=&\lim_{\mu\rightarrow +\infty}|\beta_{jk}| \int_{\R^N}\big(z_{1}^i(y)\big)^{\alpha}\big(z_{\mu^{j-i}}^j(y)\big)^{\beta}dy=0,\;\;\; 1\leq i\neq j\leq r.
\endaligned\ee
 Then conditions \eqref{zwm-aug-1}, \eqref{2con2} and \eqref{2bbar} in Lemma \ref{2lm2} hold for $(z^{1}_{\mu^1},...,z^{r}_{\mu^r})$ when $\mu>0$
 is sufficiently large. Therefore, there exists some positive constants  $\{t_j(\mu)\}_{j=1}^r$ such that $(t_1(\mu)z^{1}_{\mu^1},...,t_r(\mu)z^{r}_{\mu^r})\in \mathcal{N}$. By Lemma \ref{2lm2}, and in view of \eqref{2hardyeq} (hence ${\|z_{\mu^i}^i\|_{\la_i}^2}={|z_{\mu^i}^i|^{2^*}_{2^*}})$, we get that $t_j(\mu)\rightarrow 1$ as $\mu\rightarrow \infty, j=1,...,r.$ By \eqref{2JJJ}  and \eqref{zwm=111}, we see that
 \be   \Theta\leq J(t_1(\mu)z^1_{\mu^1},...,t_r(\mu)z^r_{\mu^r})=\frac{1}{N} \sum_{j=1}^rt^2_j(\mu)\|z^{j}_{\mu^j}\|_{\la_j}^2
 =\frac{1}{N}  \sum_{j=1}^rt^2_j(\mu) S(\la_j)^{N/2}. \nonumber
 \ee
 Letting $\mu\rightarrow \infty$ in the above equation, we get that
 \be \Theta\leq \frac{1}{N}  \sum_{j=1}^r S(\la_j)^{N/2}.  \nonumber
 \ee
On the other hand, for any $(u_1,...,u_r)\in \mathcal{N}$, we see from $\beta_{jk}<0 (j\neq k)$ and \eqref{zwm=110} that
\be \|u\|_{\la_j}^2\leq \int_{\R^N}|u|^{2^*}  \leq S(\la_j)^{-2^*/2}\|u\|_{\la_j}^{2^\ast}, \;\; j=1,...,r.  \nonumber
\ee
Combining these with \eqref{zwm-aug-15} and \eqref{2JJJ}, we get that
\be   \Theta \geq \frac{1}{N}  \sum_{j=1}^r S(\la_j)^{N/2}.  \nonumber
\ee
Hence
\be\label{zwm-aug-theta}   \Theta= \frac{1}{N}  \sum_{j=1}^r S(\la_j)^{N/2}.
\ee
Now we assume that $\Theta$ is attained by some $(u_1,...,u_r)\in \mathcal{N}$, then $(|u_1|,...,|u_r|)\in \mathcal{N}$
and $J(|u_1|,...,|u_r|)=\Theta$. By Lemma \ref{2lm1}, we know that $(|u_1|,...,|u_r|)$ is a nontrivial solution of
\eqref{model}. By the maximum principle, we may assume that $u_j>0$  in $\R^N\setminus \{0\}$ for all  $j=1,...,r$. It follows that
\be \int_{\R^N}u_i^{\alpha} u_j^{\beta}>0\;\;\hbox{for any}\;\;1\leq i,j\leq r. \nonumber
\ee
Then
\be  \|u\|_{\la_j}^2< \int_{\R^N}|u|^{2^*}dx\leq S(\la_j)^{-2^*/2}\|u\|_{\la_j}^{2^\ast}, \;\;  j=1,...,r.  \nonumber
\ee
Therefore,  it is easy to see that
$$\Theta=J(u_1,...,u_r) >\frac{1}{N}  \sum_{j=1}^r S(\la_j)^{N/2},
$$
which contradicts with \eqref{zwm-aug-theta}. This completes the proof of Theorem \ref{2th1}-(1).


\vskip0.30in

\noindent{\bf The proof of Theorem \ref{2th1}-(2)}. Note  the assumption $\beta_{jk}>0,\;  \forall j\neq k$. In this part, we  define
\be  \Theta':=\inf_{(u_1,...,u_r)\in \mathcal{N}'} J(u_1,...,u_r), \nonumber
\ee
where
\be  \mathcal{N}':=\{(u_1,...,u_r)\in \mathbb{D}\setminus \{(0,...,0)\}:   \big\langle J'(u_1,...,u_r), (u_1,...,u_r)\big\rangle=0\}.
\ee
Note that $\mathcal{N}\subset \mathcal{N}' $ and then  $\Theta'\leq \Theta$. It is easy to prove that $\Theta'>0$. Moreover, it is standard to prove that
\begin{equation}
 \aligned
  \Theta'&=\inf_{(u_1,...,u_r)\in \mathbb{D}\setminus \{(0,...,0)\}} \max_{t>0} J(tu_1,...,tu_r)\\
&=\inf_{(u_1,...,u_r)\in \mathbb{D}\setminus \{(0,...,0)\}} \frac{1}{N}\Big[\frac{  \int_{\R^N} E(u_1,...,u_r)}{ \int_{\R^N} F(u_1,...,u_r)}\Big]^{\frac{N}{2}},
\endaligned\end{equation}
here we denote that
\be E(u_1,...,u_r):=\sum_{j=1}^r \Big(|\nabla u_j|^2-\frac{\la_j}{|x|^2}u_j^2\Big), \nonumber
\ee
\be F(u_1,...,u_r):=\sum_{j=1}^r  |u_j|^{2^*}+\sum_{1\leq j<k \leq r}   2^* \beta_{jk}|u_j|^{\alpha_{jk}}|u_k|^{\alpha_{kj}}. \nonumber
\ee
Then
\be\label{zwm-aug-999}    \int_{\R^N} E(u_1,...,u_r)\geq (N \Theta')^{\frac{2}{N}}\Big(   \int_{\R^N} F(u_1,...,u_r) \Big)^{\frac{2}{2^*}},\;\;  \forall (u_1,...,u_r) \in \mathbb{D}.
\ee
Before  continuing to prove    Theorem \ref{2th1}-$(2)$, we have to establish  three lemmas. The following lemma is the counterpart of the Brezis-Lieb Lemma (see \cite{Brezis1983a})(see also  \cite{Chen2015}), here we omit the  proof.

\vskip0.2in

\bl\label{2lemma31} Let $\Omega\subset \R^N$  be an open set and $(u_n,v_n)$ be a bounded sequence in $L^{2^*}(\Omega)\times L^{2^*}(\Omega)$ .
If $(u_n, v_n)\rightarrow (u,v)$ almost everywhere in $\Omega$, then
\be\lim_{n\rightarrow \infty} \int_{\Omega}(|u_n|^{\alpha}|v_n|^{\beta}-|u_n-u|^{\alpha}|v_n-v|^{\beta}) =\int_{\Omega}|u|^{\alpha}|v|^{\beta}, \ee
here $\alpha+\beta=2^*, \alpha>1, \beta>1$.
\el

\vskip0.2in

The following lemma is the counterpart of Lions'  concentration-compactness
principle (\cite{Lions1985,Lions1985a}) for the system \eqref{model}.

\vskip0.2in

\bl\label{2 th4 lemma}   Let $(u_1,...,u_n)\in \mathbb{D}$  be a sequence such that
\be\label{2th4 0} \begin{cases}
(u_{1,n},u_{2,n},...,u_{r,n})\rightharpoonup (u_1,u_2,...,u_r) \;\;\hbox{ weakly}\; \hbox{in}\;\;\mathbb{D},\\
(u_{1,n},u_{2,n},...,u_{r,n})\rightarrow (u_1,u_2,...,u_r)\;\; \hbox{almost\ everywhere\ in}\;\; \R^N ,\\
E(u_{1,n}-u_1,u_{2,n}-u_2,...,u_{r,n}-u_r)\rightharpoonup \mu \;\;\hbox{in\ the\ sense\ of\ measures},\\
F(u_{1,n}-u_1,u_{2,n}-u_2,...,u_{r,n}-u_r)\rightharpoonup \rho\;\;\hbox{  in\ the\ sense \ of\ measures}.
 \end{cases}
 \ee
Define
\be \mu_{\infty}:=\lim_{R\rightarrow\infty}\limsup_{n\rightarrow \infty} \int_{|x|\geq R}E(u_{1,n},u_{2,n},...,u_{r,n})dx,
\ee
\be \rho_{\infty}:=\lim_{R\rightarrow\infty}\limsup_{n\rightarrow \infty} \int_{|x|\geq R}F(u_{1,n},u_{2,n},...,u_{r,n})dx.
\ee
Then it follows that
\be\label{2th4 1}\|\mu\| \geq (N\Theta')^{\frac{2}{N}} \|\rho\|^{\frac{2}{2^*}},\ee
\be\label{2th4 2}\mu_{\infty} \geq (N\Theta')^{\frac{2}{N}}{\rho_{\infty}}^{\frac{2}{2^*}}, \ee
\be\label{2th4 3}\limsup_{n\to\infty} \int_{\R^N}E(u_{1,n},...,u_{r,n})dx=\int_{\R^N}E(u_1,...,u_r)dx +\|\mu\|+\mu_{\infty},\ee
\be\label{2th4 4}\limsup_{n\rightarrow \infty} \int_{\R^N}F(u_{1,n},...,u_{r,n})dx=\int_{\R^N}F(u_1,...,u_r)dx +\|\rho\|+\rho_{\infty}. \ee
Moreover, if $(u_1,u_2,...,u_r)=(0,0,...,0)$   and $\|\mu\|=(N\Theta')^{\frac{2}{N}} \|\rho\|^{\frac{2}{2^*}}$ ,
then $\mu$ and $\rho$  are concentrated at a single point.
 \el

 \bp  In this proof we mainly follow the argument of \cite{Willem1996, Chen2015}. Firstly we assume that $(u_1,u_2,...,u_n)=(0,0,...,0)$. For any
 $h\in C_{0}^{\infty}(\R^N)$ , we see from \eqref{zwm-aug-999} that
 \be\label{2lemma4ineq1}  \int_{\R^N}E(h u_{1,n}, h u_{2,n} ,..., h u_{r,n})dx\geq (N\Theta')^{\frac{2}{N}}\Big( \int_{\R^N} |h|^{2^*}F(u_{1,n}, u_{2,n},...,u_{r,n})dx\Big)^{\frac{2}{2^*}}.
 \ee
Since $u_{j,n}\rightarrow 0,j=1,...,n$ in $L_{loc}^2(\R^N)$, we have that
\be  \int_{\R^N}E(h u_{1,n},h u_{2,n},...,h u_{r,n})dx-\int_{\R^N}|h|^2 E(u_{1,n},u_{2,n},...,u_{r,n})dx\rightarrow 0, n\rightarrow \infty.
\ee
Then by letting $n\rightarrow \infty$ in \eqref{2lemma4ineq1}, we obtain
\be\label{2 2.18}  \int_{\R^N}|h|^2 d\mu \geq (N\Theta')^{\frac{2}{N}}\Big(\int_{\R^N} |h|^{2^*}d\rho\Big)^{\frac{2}{2^*}},
\ee
that is, \eqref{2th4 1} holds.  For $R>1$, let $\psi_{R}\in C^1(\R^N)$ be such that $0\leq \psi_{R}\leq1,\psi_{R}=1$ for $|x|\geq R+1$
and $\psi_{R}=0$   for $|x|\leq R$. Then we see from \eqref{2lemma4ineq1}  that
\begin{equation*}
 \aligned & \int_{\R^N}E(\psi_{R} u_{1,n}, \psi_{R} u_{2,n},...,\psi_{R} u_{r,n})dx\\
 &\geq  (N\Theta')^{\frac{2}{N}}\Big(\int_{\R^N}  |\psi_{R}|^{2^*}F(u_{1,n}, u_{2,n},..., u_{r,n})dx\Big)^{\frac{2}{2^*}}.
\endaligned\end{equation*}
Since $u_{j,n}\rightarrow 0$  in $L_{loc}^2(\R^N)$  as $n\to \infty$   for all $j=1,...,r$, then
 \begin{equation}\label{2 2.19}
 \aligned
 & \limsup_{n\rightarrow \infty}  \int_{\R^N}|\psi_{R}|^2 E(u_{1,n},u_{2,n},...,u_{r,n})dx\\
 &\geq(N \Theta')^{\frac{2}{N}}\limsup_{n\rightarrow\infty}\Big(\int_{\R^N}|\psi_{R}|^{2^*} F(u_{1,n},u_{2,n},...,u_{r,n}) dx\Big)^{\frac{2}{2^*}}.
 \endaligned\end{equation}
Note that
\be \aligned  \int_{|x|\geq R+1}F(u_{1,n}, u_{2,n},...,u_{r,n})dx &\leq\int_{\R^N}|\psi_R|^{2^*}F(u_{1,n},u_{2,n},...,u_{r,n})dx\\
&\leq \int_{|x|\geq R}F(u_{1,n},u_{2,n},...,u_{r,n})dx,
\endaligned\ee
so
\be\label{2 2.20}  \rho_{\infty}=\lim_{R\rightarrow \infty}\limsup_{n\rightarrow \infty}\int_{\R^N}|\psi_R|^{2^*}F(u_{1,n},u_{2,n},...,u_{r,n})dx.
\ee
On the other hand,
\be \aligned
&\limsup_{n\rightarrow \infty}\int_{\R^N}|\psi_R|^2 E(u_{1,n},u_{2,n},...,u_{r,n})dx\\
&=\limsup_{n\rightarrow \infty}\int_{|x|\geq R+1}E(u_{1,n},u_{2,n},...,u_{r,n})dx\\
&\quad +\limsup_{n\rightarrow \infty}\int_{R\leq |x|\leq R+1}|\psi_R|^2 E(u_{1,n},u_{2,n},...,u_{r,n})dx\\
&=\limsup_{n\rightarrow \infty}\int_{|x|\geq R+1}E(u_{1,n},u_{2,n},...,u_{r,n})dx\\
&\quad +\limsup_{n\rightarrow \infty}\int_{R\leq |x|\leq R+1}|\psi_R|^2\Big(\sum_{j=1}^r|\nabla u_{j,n}|^2\Big)dx\\
& \geq  \limsup_{n\rightarrow \infty}\int_{|x|\geq R+1}E(u_{1,n},u_{2,n},...,u_{r,n})dx.
\endaligned\ee
 Letting $R\rightarrow \infty $ in the above inequality, we have that
\be \mu_{\infty}\leq \lim_{R\rightarrow \infty} \limsup_{n\rightarrow\infty}\int_{\R^N}
|\psi_R|^2 E(u_{1,n},u_{2,n},...,u_{r,n})dx.
\ee
Similarly,
\be \aligned
&\limsup_{n\rightarrow\infty}\int_{\R^N}|\psi_R|^2 E(u_{1,n},u_{2,n},...,u_{r,n})dx\\
& = \limsup_{n\rightarrow\infty}\int_{|x|\geq R} E(u_{1,n},u_{2,n},...,u_{r,n})dx\\
&\quad -\liminf_{n\rightarrow \infty}\int_{R\leq |x|\leq R+1}
(1-|\psi_R|^2)E(u_{1,n},u_{2,n},...,u_{r,n})dx\\
& \leq  \limsup_{n\rightarrow\infty}\int_{|x|\geq R} E(u_{1,n},u_{2,n},...,u_{r,n})dx.
 \endaligned\ee
Letting $R\rightarrow \infty$, we see that
\be \mu_{\infty}\geq \lim_{R\rightarrow \infty} \limsup_{n\rightarrow\infty}\int_{\R^N}
|\psi_R|^2 E(u_{1,n},u_{2,n},...,u_{r,n})dx.
\ee
Hence,
\be\label{2 2.21} \mu_{\infty}= \lim_{R\rightarrow \infty} \limsup_{n\rightarrow\infty}\int_{\R^N}
|\psi_R|^2 E(u_{1,n},u_{2,n},...,u_{r,n})dx.
\ee
Then \eqref{2th4 2} follows directly from \eqref{2 2.19}, \eqref{2 2.20} and \eqref{2 2.21}.
Assume moreover that $\|\mu\|=(N\Theta')^{\frac{2}{N}}\|\rho\|^{\frac{2}{2^*}}$, then by the  H\"older inequality and \eqref{2 2.18}, we have
\be\int_{\R^N}|h|^{2^*}d\rho\leq(N\Theta')^{-\frac{2}{N-2}}\|\mu\|^{\frac{2}{N-2}}\int_{\R^N}|h|^{2^*}d\mu,\;\;h\in C_0^{\infty}(\R^N).
\ee
From this we deduce that
\be\rho=(N\Theta')^{-\frac{2}{N-2}}\|\mu\|^{\frac{2}{N-2}} \mu.
\ee
So $\mu=(N\Theta')^{\frac{2}{N}}\|\rho\|^{-\frac{2}{N}}\rho$, and we see from \eqref{2 2.18} that
\be \|\rho\|^{\frac{2}{N}}(\int_{\R^N}|h|^{2^*}d\rho)^{\frac{2}{2^*}}\leq \int_{\R^N}|h|^2d\rho, \quad h\in C_0^{\infty}(\R^N).
\ee
That is, for each open set $\Omega$, we have $\rho(\Omega)^{\frac{2}{2^*}}\rho(\R^N)^{\frac{2}{N}}\leq \rho(\Omega)$. Therefore, $\rho$
is concentrated at a single point.

\vskip0.12in

For the general case, we denote that $\omega_{j,n}=u_{j,n}-u_j,j=1,2,...,r$. Then
$(\omega_{1,n}, \omega_{2,n},..., \omega_{r,n})\rightharpoonup (0,0,...,0)$ weakly in $\mathbb{D}$. From the Brezis-Lieb Lemma,
for any nonnegative function $h\in C_0(\R^N)$, we obtain  that
\be\aligned \int_{\R^N} h E(u_{1},u_{2},...,u_{r})dx=&\lim_{n \rightarrow \infty}\Big(\int_{\R^N} h E(u_{1,n},u_{2,n},...,u_{r,n})dx\\
&-\int_{\R^N} h E(\omega_{1,n},\omega_{2,n},...,\omega_{r,n})dx\Big),
\endaligned\ee
\be\aligned \int_{\R^N} h F(u_{1},u_{2},...,u_{r})dx=&\lim_{n \rightarrow \infty}\Big(\int_{\R^N} h F(u_{1,n},u_{2,n},...,u_{r,n})dx\\
&-\int_{\R^N} h F(\omega_{1,n},\omega_{2,n},...,\omega_{r,n})dx\Big),
\endaligned\ee
it follows that
\be\aligned\label{2 2.22} E(u_{1,n},u_{2,n},...,u_{r,n})&\rightharpoonup E(u_{1},u_{2},...,u_{r})+\mu,\\
F(u_{1,n},u_{2,n},...,u_{r,n})&\rightharpoonup F(u_{1},u_{2},...,u_{r})+\rho,\\
\endaligned\ee
in the sense of measures. Inequality \eqref{2th4 1} follows from the corresponding one for $(w_{1,n},w_{2,n},...,w_{r,n}).$
From the Brezis-Lieb Lemma, it is easy to prove that
\be \mu_{\infty}:=\lim_{R\rightarrow \infty}\limsup_{n\rightarrow \infty}\int_{|x|\geq R} E(w_{1,n}, w_{2,n},..., w_{r,n})dx,
\ee
\be \rho_{\infty}:=\lim_{R\rightarrow \infty}\limsup_{n\rightarrow \infty}\int_{|x|\geq R} F(w_{1,n}, w_{2,n},..., w_{r,n})dx.
\ee
Then the inequality \eqref{2th4 2} can be proved in a similar way.
 For any $R>1$, we deduce from \eqref{2 2.22} that
 \begin{equation*}
 \aligned
  &\limsup_{n\rightarrow \infty}\int_{\R^N}F(u_{1,n},u_{2,n},...,u_{r,n})\\
  &=\limsup_{n\rightarrow \infty}\Big(\int_{\R^N}|\psi_{R}|^{2^*}F(u_{1,n},u_{2,n},...,u_{r,n})\\
  & \quad +\int_{\R^N}(1-|\psi_{R}|^{2^*})F(u_{1,n},u_{2,n},...,u_{r,n})\Big)\\
  &=\limsup_{n\rightarrow \infty}\int_{\R^N}|\psi_{R}|^{2^*}F(u_{1,n},u_{2,n},...,u_{r,n})
  +\int_{\R^N}(1-|\psi_{R}|^{2^*})F(u_{1},u_{2},...,u_{r})\\
&\quad +\int_{\R^N}(1-|\psi_{R}|^{2^*})d\rho.
  \endaligned\end{equation*}
Letting $R\rightarrow \infty$, we see from \eqref{2 2.20} that \eqref{2th4 4} hold. The proof of \eqref{2th4 3} is similar.
This completes the proof of Lemma \ref{2 th4 lemma}.\ep


\vskip0.23in

\bl\label{2lemma5} Let $\beta_{jk}>0$, for any   $  1\leq j\neq k\leq r$, then \eqref{model} has a solution $(u_1, u_2,..., u_r)\in \mathbb{D}\setminus\{(0,0,...,0)\}$
(possibly  semi-trivial) such that $J(u_1,u_2,...,u_r)=\Theta'$ and that $u_j\geq 0 (j=1,...,r)$ are radially symmetric with respect to the origin.
Moreover, if further
\be\label{zwm-aug-a2}\Theta'<\min_{l\in\{1,2,...,r\}}\min_{j\in A^l}B_{j,l}^{-\frac{N-2}{2}}\frac{1}{N}S(\la_j)^{\frac{N}{2}},\ee
then $(u_1,u_2,...,u_r)\in \mathbb{D}$ is a positive ground state
solution of \eqref{model}  and $$\Theta'=J(u_1,u_2,...,u_r), \hbox{  where } B_{j,l}=\sum_{k\neq j,k\in A^l}\beta_{jk}\alpha_{jk}+1, A^l=\{1,2,...,r\}\setminus\{l\}.$$
\el

\bp For $(u_1,u_2,...,u_r)\in \mathcal{N}'$ with $u_j\geq0, j=1,...,r$, we denote by $(u_1^*,u_2^*,...,u_r^*)$ for its Schwartz symmetrization.
Then by the properties of Schwartz symmetrization, we see from $\la_j>0,\beta_{jk}>0,j\neq k$ that
\be \sum_{j=1}^r\int_{\R^N}|\nabla u_j^*|^2-\frac{\la_j}{|x|^2}|u_j^*|^2
\leq \sum_{j=1}^r\int_{\R^N}|u_j^*|^{2^*}+\sum_{1\leq j< k\leq r}2^*\beta_{jk}|u_j^*|^{\alpha_{jk}}|u_k^*|^{\alpha_{kj}}. \nonumber
\ee
Therefore, there exists $0<t^*\leq1$ such that $(t^*u_1^*, t^*u_2^*,..., t^*u_r^*)\in \mathcal{N}'$, and that
 \begin{equation}\label{2 2.23}
 \aligned
  J(t^*u_1^*,t^*u_2^*,...,t^*u_r^*)&=\frac{1}{N}(t^*)^2(\sum_{j=1}^r\|u_j^*\|_{\la_j}^2)\\
  &\leq \frac{1}{N}(\sum_{j=1}^r\|u_j^*\|_{\la_j}^2)=J(u_1,u_2,...,u_r).
 \endaligned\end{equation}
We can take a minimizing sequence $(\widetilde{u}_{1,n}, \widetilde{u}_{2,n},..., \widetilde{u}_{r,n})\in \mathcal{N}'$ such that
\be(\widetilde{u}_{1,n}, \widetilde{u}_{2,n},..., \widetilde{u}_{r,n})=(\widetilde{u}_{1,n}^*, \widetilde{u}_{2,n}^*,..., \widetilde{u}_{r,n}^*),\ee
and $J(\widetilde{u}_{1,n}, \widetilde{u}_{2,n},..., \widetilde{u}_{r,n})\rightarrow \Theta'$ as $n\rightarrow \infty$.
Define the Levy concentration function
\be Q_n(R):=\sup_{y\in\R^N}\int_{B(y,R)}F(\widetilde{u}_{1,n}, \widetilde{u}_{2,n},..., \widetilde{u}_{r,n})dx.  \nonumber
\ee
Since $\widetilde{u}_{j,n}\geq0 \;(j=1,2,...,r)$ are radially   nonincreasing, we have that
\be Q_n(R_n)=\int_{B(0,R)}F(\widetilde{u}_{1,n}, \widetilde{u}_{2,n},..., \widetilde{u}_{r,n})dx. \nonumber\ee
Then there exists $R_n>0$ such that
\be Q_n(R)=\int_{B(0,R_n)}F(\widetilde{u}_{1,n}, \widetilde{u}_{2,n},..., \widetilde{u}_{r,n})dx
=\frac{1}{2}\int_{\R^N}F(\widetilde{u}_{1,n}, \widetilde{u}_{2,n},..., \widetilde{u}_{r,n})dx.
\nonumber \ee
Define
\be\aligned &\big(u_{1,n}(x),u_{2,n}(x),...,u_{r,n}(x)\big):\nonumber\\
&=\big(R_n^{\frac{N-2}{2}}\widetilde{u}_{1,n}(R_n x),R_n^{\frac{N-2}{2}}\widetilde{u}_{2,n}(R_n x),...,R_n^{\frac{N-2}{2}}\widetilde{u}_{r,n}(R_n x)\big).
\endaligned
\nonumber\ee
By a direct computation, we know that $(u_{1,n}, u_{2,n},..., u_{r,n})\in \mathcal{N}', J(u_{1,n}, u_{2,n},..., u_{r,n})\rightarrow \Theta'$ and that  $ u_j\geq 0$
are radially nonincreasing. Moreover,
\be\label{zwm-aug-z7} \aligned \int_{B(0,1)}F(u_{1,n}, u_{2,n},..., u_{r,n})dx&=\frac{1}{2}\int_{\R^N}F(u_{1,n}, u_{2,n},...,u_{r,n})dx \\
&=\sup_{y\in\R^N}\int_{B(y,1)}F(u_{1,n}, u_{2,n},..., u_{r,n})dx.
\endaligned
 \ee
From \eqref{2 2.23}, we know that $(u_{1,n}, u_{2,n},..., u_{r,n})$ are uniformly bounded in $\mathbb{D}$. Then passing to a
subsequence, there exist $(u_1,u_2,...,u_r)\in \mathbb{D}$ and finite measures $\mu,\rho$ such that \eqref{2th4 0} holds.
Then by Lemma \ref{2 th4 lemma}, we see that \eqref{2th4 1}-\eqref{2th4 4} hold. Note that
\be\label{2 2.25}\sum_{j=1}^r \|u_{j,n}\|_{\la_j}^2=\int_{\R^N} F(u_{1,n}, u_{2,n},..., u_{r,n})dx\rightarrow N\Theta',\;\;\hbox{as}\ n\rightarrow\infty.
\nonumber \ee
From \eqref{2th4 0}-\eqref{2th4 4}, we have that
\be\label{zwm-aug-z8} N\Theta'=\int_{\R^N} F(u_1,u_2,...,u_r)dx+\|\rho\|+\rho_{\infty},
\ee
\be\label{zwm-aug-z9} N\Theta'\geq(N\Theta')^{\frac{2}{N}}[(\int_{\R^N} F(u_1,u_2,...,u_r)dx)^{\frac{2}{2^*}}+\|\rho\|^{\frac{2}{2^*}}+\rho_{\infty}^{\frac{2}{2^*}}].
\ee
Therefore, $\int_{\R^N} F(u_1,u_2,...,u_r)dx$, $\rho$ and $\rho_{\infty}$ are equal to either 0 or $N\Theta'$. By \eqref{zwm-aug-z7}, \eqref{zwm-aug-z8}-\eqref{zwm-aug-z9},
we have $\rho_{\infty}\leq \frac{1}{2}N\Theta'$, hence, $\rho_{\infty}=0$. If $\|\rho\|=N\Theta'$, then $$\int_{\R^N} F(u_1,u_2,...,u_r)dx=0$$ and so
$(u_1,u_2,...,u_r)=(0,0,...,0)$. On the other hand, since $\|\mu\|\leq N\Theta'$, we deduce from \eqref{2th4 1}
that $\|\mu\|=(N\Theta')^{\frac{2}{N}}\|\rho\|^{\frac{2}{2^*}}$. Then Lemma \ref{2 th4 lemma} implies that $\rho$ is
concentrated at a single point $y_0$, and we see from \eqref{zwm-aug-z7}, \eqref{zwm-aug-z8} and \eqref{zwm-aug-z9}, that
\be\aligned \frac{1}{2}N\Theta'&=\lim_{n \rightarrow \infty}\sup_{y\in\R^N}\int_{B(y,1)}F(u_{1,n},u_{2,n},...,u_{r,n})\\
&\geq\lim_{n\rightarrow \infty}\int_{B(y_0,1)}F(u_{1,n}, u_{2,n},..., u_{r,n})=\|\rho\|,
\endaligned
\nonumber\ee
a contradiction. Therefore, $\int_{\R^N}F(u_1,u_2,...,u_r)dx=N\Theta'$. Since $\sum_{j=1}^r\|u_j\|_{\la_j}^2\leq N\Theta'$, we deduce from
\eqref{zwm-aug-z7}, \eqref{zwm-aug-z8} and \eqref{zwm-aug-z9}  that
\be N\Theta'=\sum_{j=1}^r\|u_j\|_{\la_j}^2=\int_{\R^N}F(u_1,u_2,...,u_r)dx=\lim_{n\rightarrow \infty}\sum_{j=1}^r\|u_{j,n}\|_{\la_j}^2,
\nonumber \ee
then
\be \sum_{j=1}^r\|u_{j,n}-u_j\|_{\la_j}^2=\sum_{j=1}^r\|u_{j,n}\|_{\la_j}^2+\|u_j\|_{\la_j}^2
-2\langle u_{j,n},u_j\rangle_{\la_j}\rightarrow 0,\;\hbox{as} \ n\rightarrow\infty, \nonumber
\ee
that is, $(u_{1,n}, u_{2,n},..., u_{r,n})\rightarrow (u_1, u_2,..., u_r)$ strongly in $\mathbb{D}$, $(u_1,u_2,...,u_r)\in \mathcal{N}'$
and $J(u_1, u_2,..., u_r)=\Theta'$.
 Recall that $\Theta'>0$, hence $(u_1, u_2,..., u_r)\neq (0, 0 ,..., 0)$. By the definition of $\mathcal{N}'$
and using the Lagrange multiplier method, it is standard to prove that $J(u_1, u_2,..., u_r)=0$. Therefore, $(u_1, u_2 ,..., u_r)$ is a solution of
\eqref{model}.

\vskip0.3in

Now assume that $\Theta'<\min_{l\in\{1,2,...,r\}}\min_{j\in A^l}B_{j,l}^{-\frac{N-2}{2}}\frac{1}{N}S(\la_j)^{\frac{N}{2}}$.
Then it is easy to prove that
\be u_j\not\equiv 0,\;\;  j=1,2,...,r. \nonumber
\ee
In fact,  note that  $\alpha_{jk}+\alpha_{kj}=2^*,$
\be \int_{\R^N}|u_j|^{\alpha_{jk}}|u_k|^{\alpha_{kj}}\leq \frac{\alpha_{jk}}{2^*}\int_{\R^N}|u_j|^{2^*}+
\frac{\alpha_{kj}}{2^*}\int_{\R^N}|u_k|^{2^*}, \nonumber
\ee
If $u_{l}\equiv0$, then
\be \aligned
& \sum_{j=1}^r\int_{\R^N}|u_j|^{2^*}+\sum_{1\leq j<k\leq r}2^*\beta_{jk}|u_j|^{\alpha_{jk}}|u_k|^{\alpha_{kj}}\\
&\leq \sum_{j\in A^l}\Big(\sum_{k\neq j,k\in A^l}\beta_{jk}\alpha_{jk}+1\Big)\int_{\R^N}|u_j|^{2^*}\\
&\leq \sum_{j\in A^l} B_{j,l}\int_{\R^N}|u_j|^{2^*},
\endaligned\ee
where $A^l=\{1, 2,..., r\}\setminus\{l\}$,
then
\be\label{2I'}\aligned \Theta'&=\inf_{(u_1,u_2,...,u_j)\in \mathbb{D}\setminus \{(0, 0,..., 0)\}}\\
& \quad\quad \frac{1}{N}
\Big[\frac{\sum_{j=1}^r\int_{\R^N}\|u_j\|_{\la_j}^2}{(\sum_{j=1}^r\int_{\R^N}|u_j|^{2^*}+\sum_{1\leq j<k\leq r}2^*\beta_{jk}|u_j|^{\alpha_{jk}}|u_k|^{\alpha_{kj}})^{\frac{2}{2^*}}}\Big]^{\frac{N}{2}}\\
&\geq \inf_{(u_1, u_2, ..., u_j)\in \mathbb{D}\setminus \{(0, 0, ..., 0)\}}\frac{1}{N}
\Big[\frac{\sum_{j\in A^l\|u_j\|_{\la_j}^2}}{(\sum_{j\in A^l} B_{j,l}\int_{\R^N}|u_j|^{2^*})^{\frac{2}{2^*}}}\Big]^{\frac{N}{2}}\\
&\geq \min_{j\in A^l}B_{j,l}^{-\frac{N-2}{2}}\frac{1}{N}S(\la_j)^{\frac{N}{2}}\\
&\geq \min_{l\in\{1, 2, ..., r\}}\min_{j\in A^l}B_{j, l}^{-\frac{N-2}{2}}\frac{1}{N}S(\la_j)^{\frac{N}{2}},
\endaligned\ee
a contradiction. Hence,  $u_j\not\equiv 0$ for all $ j=1,2,...,r$. That is,
 $(u_1,u_2,...,u_r)\in \mathcal{N}'$ and so $J(u_1,u_2,...,u_r)=\Theta'=\Theta$. Hence,
$(u_1,u_2,...,u_r)$ is a ground state solution of \eqref{model}. By the maximum principle, $u_j>0$ in
$\R^N\setminus \{0\}$. This completes the proof of Lemma \ref{2lemma5}.
\ep

\vskip0.23in
\noindent {\bf Finishing the proof of Theorem \ref{2th1}-(2)}.  We apply Lemma \ref{2lemma5}.  It suffices to prove \eqref{zwm-aug-a2}.
Let $\beta_{jk}>0,j\neq k$.
Without loss of generality, we assume that $\la_1\leq \la_2\leq\cdots \leq\la_r,$ then $S(\la_r)\leq \cdots  S(\la_2)\leq S(\la_1)$.
Denote that $d_{\alpha}:=\frac{\La_N-\la_{\alpha}}{\La_N-\la_r}$, for $1\leq\alpha \leq r$. By Hardy's inequality, we know  that
$\|u\|_{\la_{\alpha}}\leq d_{\alpha}\|u\|_{\la_r}$ for all $u\in D^{1,2}(\R^N), \alpha=1,2,...,r.$ Then we deduce from
\eqref{2I'} that
 \begin{equation}
 \aligned
\Theta'&\leq \frac{1}{N}\Big[\frac{\sum_{j=1}^r \|z^r_{\mu}\|^2_{\la_j}}{(\sum_{j=1}^r|z^r_{\mu}|^{2^*}
+\sum_{1\leq k<j\leq r}2^*\beta_{jk}|z^r_{\mu}|^{2^*})^{\frac{2}{2^*}}}\Big]^{\frac{N}{2}}\\
&\leq \frac{1}{N}\Big[\frac{(1+\sum_{\alpha=1}^{r-1}d_{\alpha})}{(r+\sum_{j,k=1,j\neq k}^r\frac{2^*}{2}\beta_{jk})^{\frac{2}{2^*}}}
\frac{\|z^r_{\mu}\|^2_{\la_r}}{(\int_{\R^N}|z^r_{\mu}|^{2^*})^{\frac{2}{2^*}}}\Big]^{\frac{N}{2}}\\
&=\Big[\frac{(1+\sum_{\alpha=1}^{r-1}d_{\alpha})}{(r+\sum_{j,k=1,j\neq k}^r\frac{2^*}{2}\beta_{jk})^{\frac{2}{2^*}}}
\Big]^{\frac{N}{2}}\frac{1}{N}S(\la_r)^{\frac{N}{2}}\\
&<\min_{l\in\{1,2,...,r\}}\min_{j\in A^l}B_{j,l}^{-\frac{N-2}{2}}\frac{1}{N}S(\la_j)^{\frac{N}{2}}
 \endaligned\end{equation}
 provided that  $$(r+\sum_{j,k=1,j\neq k}^r\frac{2^*}{2}\beta_{jk})/(\max_{j,l}B_{j,l})>(1+\sum_{\alpha=1}^{r-1}d_{\alpha})^{\frac{N}{N-2}},$$
  where,
$B_{j,l}=\sum_{k\neq j,k\in A^l}\beta_{jk}\alpha_{jk}+1$, $A^l=\{1,2,...,r\}\setminus\{l\}$.
Hence, the conclusion follows from Lemma \ref{2lemma5}.  \hfill$\Box$

\section{Proof of Theorem \ref{2th2}}

  In this section, we will use the moving plane  method to prove Theorem \ref{2th2}. In the sequel,
we assume that $ N=3$  or $N=4, r\geq 3,\alpha_{jk}+\alpha_{kj}=2^*,\alpha_{jk}\geq 2,\alpha_{kj}\geq 2,\beta_{jk}>0$ for $k\neq j$
 and $\la_j\in(0,\Lambda_N)$ for all $ j=1,...,r$.
Let $(u_1,...,u_r)$ be any   positive solution of \eqref{model}. For $\la<0$, we consider the following  reflection:
\be x=(x_1,x_2,...,x_N)\mapsto x^{\la}=(2\la-x_1, x_2,..., x_N), \nonumber \ee
where $x\in \Sigma^{\la}:=\{x\in \R^N:x_1<\la\}$. Define $u_j^{\la}(x):=u_j(x^{\la})$, then
\be u_j(x)=u_j^{\la}(x) \;\;\hbox{for}\;\;x\in \pa \Sigma^{\la},  \hbox{ where } \;\Sigma^{\la}:=\{x\in \R^N:x_1=\la\}. \nonumber \ee
Define $w_j^{\la}(x):=u_j^{\la}(x)-u_j(x) $ for $x\in \Sigma^{\la}$, then
\be w_j^{\la}=0,\;\; x\in \pa\Sigma^{\la}.  \nonumber  \ee
Recall that $u_j(x)$ and $ u_j^{ \la}(x) $ satisfy the following equations
\begin{equation}\label{2moveq}
 \aligned
-\Delta u_j-\frac{\la_j}{|x|^2}u_j&=u_j^{2^*-1}+\sum_{k\neq j}\alpha_{jk}\beta_{jk}u_j^{\alpha_{jk}-1}u_k^{\alpha_{kj}},\\
-\Delta u_j^{\la}-\frac{\la_j}{|x|^2}u_j^{\la}&=(u_j^{\la})^{2^*-1}+\sum_{k\neq j}\alpha_{jk}\beta_{jk}(u_j^{\la})^{\alpha_{jk}-1}(u_k^{\la})^{\alpha_{kj}},
 \endaligned\end{equation}
thus we have
\begin{equation}
 \aligned
 -\Delta w_j^{\la}\geq \frac{\la_j}{|x|^2}w_j^{\la}+b_{jj}^{\la}+\sum_{k \neq j} b_{jk}w_{k}^{\la}, \nonumber
 \endaligned\end{equation}
 here
\begin{equation}
 \aligned
b_{jj}^{\la}&=\frac{(u_j^{\la})^{2^*-1}-u_j^{2^*-1}}{u_j^{\la}-u_j}+
\sum_{k\neq j} \alpha_{jk}\beta_{jk}u_k^{\alpha_{kj}}\frac{(u_j^{\la})^{2^*-1}-u_j^{\alpha_{jk}-1}}{u_j^{\la}-u_j}\geq 0, \nonumber \\
b_{jk}^{\la}&=\alpha_{jk}\beta_{jk}(u_j^{\la})^{2^*-1}\frac{(u_j^{\la})^{2^*-1}-u_j^{\alpha_{jk}-1}}{u_j^{\la}-u_j}\geq 0.\nonumber
 \endaligned\end{equation}
 Define
 \be\label{zwm-aug-xw1} \Omega_j^{\la}:=\{x \in \Sigma^{\la}:w_j^{\la}(x)<0\},j=1,...,r.
 \ee
Since  $u_j\in L^{2^*}(\R^N)$ and $\Omega_j^{\la}\subset \Sigma ^{\la}$, there exists $\la_{0}\rightarrow -\infty$ such that
\begin{equation}\label{2mov0}
 \aligned
&\|b_{jj}^{\la}\|_{L^{\frac{N}{2}}(\Omega_j^{\la})}\rightarrow 0,\\
&\|b_{jk}^{\la}\|_{L^{\frac{N}{2}}(\Omega_k^{\la})}\rightarrow 0,\\
&\|u_j^{\alpha_{jk}}u_k^{\alpha_{kj}}\|_{L^{\frac{N}{2}}(\Omega_j^{\la}\bigcap\Omega_k^{\la})}\rightarrow 0.
 \endaligned\end{equation}
{\bf Step 1:}  We claim that for any $\la\leq \la_0$, all $w_j^{\la}>0$ in $\Sigma^{\la}\setminus \{0^{\la}\}.$  For this aim, we  define
\be w_{j,-}^{\la}:=-\max\{-w_j^{\la},0\},\;\;   j=1,...,r,  \nonumber \ee
then $w_{j,-}^{\la}\in D^{1,2}(\R^N)$ and
\begin{equation*}
 \aligned
\int_{\Omega_j^{\la}}|\nabla w_{j,-}^{\la}|^2&\leq\int_{\Omega_j^{\la}} \frac{\la_j}{|x|^2}| w_{j,-}^{\la}|^2+\int_{\Omega_j^{\la}} b_{jj}^{\la}|w_{j,-}^{\la}|^2+\sum_{k\neq j}\int_{\Omega_j^{\la}\bigcap\Omega_k^{\la}} b_{jk}^{\la}w_{k,-}^{\la}w_{j,-}^{\la},\\
&\leq \frac{\la_j}{\Lambda_N}\int_{\Omega_j^{\la}}|\nabla w_{j,-}^{\la}|^2+\|b_{jj}^{\la}\|_{L^{\frac{N}{2}}(\Omega_j^{\la})}\|w_{j,-}^{\la}\|^2_{L^{2^*}(\Omega_j^{\la})}\\
&\quad +\sum_{k\neq j}\|b_{jk}^{\la}\|_{L^{\frac{N}{2}}(\Omega_j^{\la}\bigcap\Omega_k^{\la})}\|w_{k,-}^{\la}\|_{L^{2^*}(\Omega_k^{\la})}
\|w_{j,-}^{\la}\|_{L^{2^*}(\Omega_j^{\la})},
 \endaligned\end{equation*}
then we have
\begin{equation*}
 \aligned
& (1-\frac{\la_j}{\Lambda_N}\int_{\Omega_j^{\la}}|\nabla w_{j,-}^{\la}|^2)\\
&\leq\|b_{jj}^{\la}\|_{L^{\frac{N}{2}}(\Omega_j^{\la})}\|w_{j,-}^{\la}\|^2_{L^{2^*}(\Omega_j^{\la})}\\
&\quad +\sum_{k\neq j}\frac{1}{2}\|b_{jk}^{\la}\|_{L^{\frac{N}{2}}(\Omega_j^{\la}\bigcap\Omega_k^{\la})}(\|w_{k,-}^{\la}\|^2_{L^{2^*}(\Omega_k^{\la})}
+\|w_{j,-}^{\la}\|^2_{L^{2^*}(\Omega_j^{\la})})\\
&=(\|b_{jj}^{\la}\|_{L^{\frac{N}{2}}(\Omega_j^{\la})}+\sum_{k\neq j}\frac{1}{2}\|b_{jk}^{\la}\|_{L^{\frac{N}{2}}(\Omega_j^{\la}\bigcap\Omega_k^{\la})})
\|w_{j,-}^{\la}\|^2_{L^{2^*}(\Omega_j^{\la})}\\
&\quad +\sum_{k\neq j}\frac{1}{2}\|b_{jk}^{\la}\|_{L^{\frac{N}{2}}(\Omega_j^{\la}\bigcap\Omega_k^{\la})}\|w_{k,-}^{\la}\|^2_{L^{2^*}(\Omega_k^{\la})}.
 \endaligned\end{equation*}
Denote that
\be a_{jj}^{\la}:=1-\frac{\la_j}{\la_N}-(\|b_{jj}^{\la}\|_{L^{\frac{N}{2}}(\Omega_j^{\la})}
+\sum_{k\neq j}\frac{1}{2}\|b_{jk}^{\la}\|_{L^{\frac{N}{2}}(\Omega_j^{\la}\bigcap\Omega_k^{\la})}),    \nonumber
\ee
then by \eqref{2mov0}, we see that by letting $|\la_0|\rightarrow\infty, \la_0<0$, we can make sure that
\be a_{jj}^{\la}\geq \frac{1}{2}(1-\frac{\la_j}{\Lambda_N})\geq \min_j \frac{1}{2} (1-\frac{\la_j}{\Lambda_N}):=a.   \nonumber
\ee
By this notation
\begin{equation}
 \aligned
a\int_{\Omega_j^{\la}}|\nabla w_{j,-}^{\la}|^2&\leq
\sum_{k\neq j}\frac{1}{2}\|b_{jk}^{\la}\|_{L^{\frac{N}{2}}(\Omega_j^{\la}\bigcap\Omega_k^{\la})}\|w_{k,-}^{\la}\|^2_{L^{2^*}(\Omega_k^j)}   \nonumber  \\
&\leq \sum_{k\neq j}\frac{1}{2}\|b_{jk}^{\la}\|_{L^{\frac{N}{2}}(\Omega_j^{\la}\bigcap\Omega_k^{\la})}
 \frac{1}{S} \int_{\Omega_k^j}|\nabla w_{j,-}^{\la}|^2.
 \endaligned\end{equation}
Summing up  the above  inequalities, we have
\begin{equation*}
 \aligned
a\sum_{j=1}^r \int_{\Omega_j^{\la}}|\nabla w_{j,-}^{\la}|^2&\leq\sum_{j=1}\sum_{k\neq j}\frac{1}{2}\|b_{jk}^{\la}\|_{L^{\frac{N}{2}}(\Omega_j^{\la}\bigcap\Omega_k^{\la})} \frac{1}{S} \int_{\Omega_k^{\la}}|\nabla w_{j,-}^{\la}|^2\\
&=\sum_{k=1}^r\sum_{j\neq k}\frac{1}{2}\|b_{jk}^{\la}\|_{L^{\frac{N}{2}}(\Omega_j^{\la}\bigcap\Omega_k^{\la})} \frac{1}{S} \int_{\Omega_k^{\la}}|\nabla w_{k,-}^{\la}|^2\\
&=\sum_{j=1}^r(\sum_{k\neq j}\frac{1}{2}\|b_{kj}^{\la}\|_{L^{\frac{N}{2}}(\Omega_j^{\la}\bigcap\Omega_k^{\la})} \frac{1}{S} )\int_{\Omega_k^{\la}}|\nabla w_{j,-}^{\la}|^2.
 \endaligned\end{equation*}
Recall  \eqref{2mov0}, we can  let $\la_0$  go to  $-\infty$ such that
\be \sum_{k\neq j}\frac{1}{2}\|b_{kj}^{\la}\|_{L^{\frac{N}{2}}(\Omega_j^{\la}\bigcap\Omega_k^{\la})} \frac{1}{S} \leq\frac{1}{2} a.   \nonumber
\ee
It follows  that
\be a\sum_{j=1}^r \int_{\Omega_j^{\la}}|\nabla w_{j,-}^{\la}|^2\leq \frac{1}{2}a\sum_{j=1}^r \int_{\Omega_j^{\la}}|\nabla w_{j,-}^{\la}|^2,  \nonumber
\ee
hence we have the following inequality
\be \frac{1}{2}a\sum_{j=1}^r \int_{\Omega_j^{\la}}|\nabla w_{j,-}^{\la}|^2\leq 0.  \nonumber
\ee
Since $a>0$, we have
\be \int_{\Omega_j^{\la}}|\nabla w_{j,-}^{\la}|^2=0, \quad  j=1,...,r.   \nonumber
\ee
Therefore,  we have the following alternative conclusion   for any $j$: either
\be\label{2moval1} w_{j,-}^{\la}(x)\equiv  const \;\; \hbox{ in }  \Omega _j^{\la},  \;\; \;\; m(\Omega _j^{\la})>0,
\ee
or
\be \label{2moval2}m(\Omega _j^{\la})=0,
\ee
where $m$ represents the  Lebesgue measure. Notice that $w_{j,-}^{\la}=0$ on $\pa \Omega_j^{\la}$, hence, by  \eqref{2moval1} we have  that
$w_{j,-}^{\la}(x)\equiv 0$ in $\Omega_j^{\la}$.  This is,  $w_j^{\la}\geq 0$ in $\Sigma^{\la}\setminus\{0^{\la}\}$.
 Look \eqref{2moval2} now,  it just says that $w_j^{\la}\geq 0$ in $\Sigma^{\la}\setminus\{0^{\la}\}$ in another way. In summary,
 we have that $w_j^{\la}\geq 0$ in $\Sigma^{\la}\setminus\{0^{\la}\}$. Recall  \eqref{2moveq},  we see that
\be -\Delta w_j^{\la}\geq \la_j(\frac{1}{|x^{\la}|^2}-\frac{1}{|x|^2}) u_j^{\la}(x)>0, \nonumber
\ee
holds in $\Sigma^{\la}\setminus\{0^{\la}\}$.  Then by the strong maximum principle, we have $w_j^{\la}> 0$ in $\Sigma^{\la}\setminus\{0^{\la}\}$.

\vskip0.1in

\noindent {\bf Step 2:}  Define $\la^*=\sup\{\overline{\la}<0:w_j^{\la}>0 $ in $\Sigma^{\la}\setminus\{0^{\la}\}, j=1,...,r, \forall \la < \overline{\la}\}$. Then we claim that $\la^*=0.$

\vskip0.11in
Assume by contradiction that $\la^*<0$. By the continuity we have $w_j^{\la^*}\geq0$ in $\Sigma^{\la}\setminus\{0^{\la}\}$.
 By a similar argument as in Step 1,   we have $w_j^{\la^*}>0, j=1,...,r$ in $\Sigma^{\la}\setminus\{0^{\la}\}$.
By the absolutely continuity of  the integral, there exists  a $\la$ with $0>\la>\la^*$ such that
\begin{equation}\label{2mov01}
 \aligned
&\|b_{jj}^{\la}\|_{L^{\frac{N}{2}}(\Omega_j^{\la})}\rightarrow 0,\\
&\|b_{jk}^{\la}\|_{L^{\frac{N}{2}}(\Omega_k^{\la})}\rightarrow 0,\\
&\|u_j^{\alpha_{jk}}u_k^{\alpha_{kj}}\|_{L^{\frac{N}{2}}(\Omega_j^{\la}\bigcap\Omega_k^{\la})}\rightarrow 0, \quad \hbox { as }   \la\rightarrow \la^*.
 \endaligned\end{equation}
Then we can follow the same proof as in Step 1, we can find a $\la$, which satisfies $0>\la>\la^*$
and $w_j^{\la}>0$ in $\Sigma^{\la}\setminus\{0^{\la}\}$ as  $\la$ closing  to $\la^*$, which contradicts to the definition of $\la^*$. Therefore, $\la^*=0$.

\vskip0.12in

\noindent {\bf Step 3.}   We show that $w_j (j=1,...,r) $ are radially symmetric with respect to the origin.
Since $\la^*=0$, then we can carry out the above procedure in the opposite direction, namely
 we can take the transform $y=(y_1,y_2,...,y_r)=(-x_1,x_2,...,x_r), $ then moving plane by
 Step 1 and Step 2 about $y_1$, we can derive that $u_j \, (j=1,...,r) $ are symmetric with respect to 0 in the $x_1$ direction.
  Since we take the orthogonal transform $y=(y_1,y_2,...,y_r)=A (x_1, x_2,..., x_r)$ arbitrarily, where $A$ is a $r$ order orthogonal matrix,
   we can derive that $u_j\, (j=1,...,r) $ are symmetric  with respect to 0 in any direction.  It follows that  $u_j\, (j=1,...,r) $
   are radially symmetric with respect to the origin. This completes the proof of Theorem   \ref{2th2}.

\section{Proofs of Theorems \ref{2th3}-\ref{2th4} }
Firstly,  we will prepare  several lemmas
which are  essential to  the  proof of Theorems \ref{2th3}-\ref{2th4}. We remark  that these lemmas are also interesting
from its own perspective.

\bl\label{2corrth4}  Consider the  following nonlinear constraint problem:

\be\label{2se3tem}\begin{cases}
x_1+x_2\leq c_1+c_2,\\
f_1(x_1,x_2):=x_1^{p-1}+\nu \alpha x_1^{\frac{\alpha}{2}-1}x_2^{\frac{\beta}{2}}\geq f_1(c_1,c_2),\\
f_2(x_1, x_2):=x_2^{p-1}+\nu \beta x_1^{\frac{\alpha}{2}}x_2^{\frac{\beta}{2}-1}\geq f_2(c_1,c_2).
\end{cases}
\ee
If  $ \nu>(p-1)/\min\{d_1(\alpha,\beta), d_2(\alpha,\beta), d_3(\alpha,\beta)\}, $
then
\be    x_1=c_1, \quad  x_2=c_2,  \nonumber \ee
where $N\geq5, p=\frac{2^*}{2}, \alpha+\beta=2^*=2p$ and $\alpha>0, \beta>0$ ,$c_i>0,x_i>0; i=1,2;$
\be d_1(\alpha,\beta)=2p(1-\frac{\alpha}{2})^{\frac{\alpha}{2p}}(1-\frac{\beta}{2})^{\frac{\beta}{2p}}\;\;\hbox{if}\;\;\alpha\neq\beta;
d_1(\alpha,\beta)=p \;\;\hbox{if}\;\;\alpha=\beta, \nonumber
\ee
\be d_2(\alpha,\beta)=\beta(1-\frac{\beta}{2})^{1-\frac{\alpha}{2}}(1-\frac{\alpha}{2})^{\frac{\alpha}{2}}
+\frac{1}{2}\alpha\beta(1-\frac{\alpha}{2})^{1-\frac{\alpha}{2}}(1-\frac{\beta}{2})^{\frac{\alpha}{2}-1}, \nonumber
\ee
\be d_3(\alpha,\beta)=\alpha(1-\frac{\alpha}{2})^{1-\frac{\beta}{2}}(1-\frac{\beta}{2})^{\frac{\beta}{2}}
+\frac{1}{2}\alpha\beta(1-\frac{\alpha}{2})^{1-\frac{\beta}{2}}(1-\frac{\beta}{2})^{\frac{\beta}{2}-1}. \nonumber
\ee
In particular if $\alpha=\beta=p$ in \eqref{2se3tem}, we have the concise form, that is,
\be  d_1(\alpha,\beta)=d_2(\alpha,\beta)=d_3(\alpha,\beta)=p,  \nonumber \ee
hence, under this case, $x_1=c_1$ and $x_2=c_2$ if $\nu>\frac{2}{N}.$
\el


\bp
By Lemma \ref{2pringe}, we only need to check that the matrix $F=(\frac{\pa f_j}{\pa x_i}):=(F_{ij})$ satisfies
$\det(F)<0, F_{22}-F_{12}<0, F_{11}-F_{21}<0.$ By a direct computation, we have
$$ F=
\begin{pmatrix}\label{2Fmatrix}
      (p-1)x_1^{p-2}+\nu\alpha(\frac{\alpha}{2}-1)x_1^{\frac{\alpha}{2}-2}x_2^{\frac{\beta}{2}}
      & \frac{1}{2}\nu\alpha\beta x_1^{\frac{\alpha}{2}-1}x_2^{\frac{\beta}{2}-1}\\
       \frac{1}{2}\nu\alpha\beta x_1^{\frac{\alpha}{2}-1}x_2^{\frac{\beta}{2}-1}
       &(p-1)x_2^{p-2}+\nu\beta(\frac{\beta}{2}-1)x_1^{\frac{\alpha}{2}}x_2^{\frac{\beta}{2}-2}

\end{pmatrix},$$

 $$ F^{-1}= \frac{1}{\det(F)}
 \begin{pmatrix}\label{2Fmatrix}
A_0
 & -\frac{1}{2}\nu\alpha\beta x_1^{\frac{\alpha}{2}-1}x_2^{\frac{\beta}{2}-1}\\
 \frac{1}{2}\nu\alpha\beta x_1^{\frac{\alpha}{2}-1}x_2^{\frac{\beta}{2}-1}
  &(p-1)x_1^{p-2}+\nu\alpha(\frac{\alpha}{2}-1)x_1^{\frac{\alpha}{2}-2}x_2^{\frac{\beta}{2}}
\end{pmatrix},$$
where $A_0:= (p-1)x_2^{p-2}+\nu\beta(\frac{\beta}{2}-1)x_1^{\frac{\alpha}{2}}x_2^{\frac{\beta}{2}-2}$ and

  \begin{equation}
  \aligned
\det(F)= &\Big\{(p-1)^2+(\nu^2\alpha\beta(\frac{\alpha}{2}-1)(\frac{\beta}{2}-1)-\frac{1}{4}\nu^2\alpha^2\beta^2)(\frac{x_1}{x_2})^{\alpha-p}\\
&+(p-1)\nu\beta(\frac{\beta}{2}-1)(\frac{x_1}{x_2})^{\frac{\alpha}{2}}
+(p-1)\nu\alpha(\frac{\alpha}{2}-1)(\frac{x_2}{x_1})^{p-\frac{\alpha}{2}}\Big\}x_1^{p-2}x_2^{p-2}.
\endaligned
\end{equation}
When  $\alpha\neq\beta$, since $(\frac{\alpha}{2}-1)(\frac{\beta}{2}-1)-\frac{1}{4}\alpha\beta=1-p<0$, we have
   \begin{equation}
  \aligned
\det(F)\leq&\Big\{(p-1)^2+(p-1)\nu\beta(\frac{\beta}{2}-1)(\frac{x_1}{x_2})^{\frac{\alpha}{2}}\\
&+(p-1)\nu\alpha(\frac{\alpha}{2}-1)(\frac{x_2}{x_1})^{p-\frac{\alpha}{2}}\Big\}x_1^{p-2}x_2^{p-2}.\\\endaligned
  \end{equation}
When $\alpha=\beta$,
   \begin{equation}
  \aligned
&\det(F)\leq\Big\{(p-1)^2+\nu^2\alpha\beta(\frac{\alpha}{2}-1)(\frac{\beta}{2}-1)-\\
&\frac{1}{4}\nu^2\alpha^2\beta^2+
(p-1)\nu\beta(\frac{\beta}{2}-1)(\frac{x_1}{x_2})^{\frac{\alpha}{2}}+(p-1)\nu\alpha(\frac{\alpha}{2}-1)
(\frac{x_2}{x_1})^{p-\frac{\alpha}{2}}\Big\}x_1^{p-2}x_2^{p-2}\\
&= (p-1)^2+\nu^2p^2(1-p)+\nu(p-1)p(\frac{p}{2}-1)(\frac{x_1}{x_2})^{\frac{p}{2}}\\
&\quad +\nu(p-1)p(\frac{p}{2}-1)(\frac{x_2}{x_1})^{\frac{p}{2}}.
  \endaligned
  \end{equation}
Let
\be h_1(x):=(p-1)\nu\beta(\frac{\beta}{2}-1)x^{\frac{\alpha}{2}}
+(p-1)\nu\alpha(\frac{\alpha}{2}-1)(\frac{1}{x})^{p-\frac{\alpha}{2}}, 0<x:=\frac{x_1}{x_2}<\infty,
\ee
hence,
\be h_1'(x_0)=0 \Rightarrow x_0=\Big[\frac{1-\frac{\alpha}{2}}{1-\frac{\beta}{2}}\Big]^{\frac{1}{p}}.
\ee
It is easy to see that $x_0$ is the maximum point of $h_1(x)$ in the interval $(0,\infty)$, so
\be\label{2cd1} \det(F)<0\Leftrightarrow h_1(x)<0\Leftrightarrow \nu>(p-1)/d_1(\alpha, \beta).
\ee
Next we estimate $F_{22}-F_{12}:$
\be  F_{22}-F_{12}=x_2^{p-2}\{(p-1)-\nu\beta(1-\frac{\beta}{2})(\frac{x_1}{x_2})^{\frac{\alpha}{2}}
-\frac{1}{2}\alpha\beta\nu(\frac{x_1}{x_2})^{\frac{\alpha}{2}-1}\}, \nonumber
\ee
\be   h_2(x):=(p-1)-\nu\beta(1-\frac{\beta}{2})x^{\frac{\alpha}{2}}
-\frac{1}{2}\alpha\beta\nu(x)^{\frac{\alpha}{2}-1}, 0<x:=\frac{x_1}{x_2}<\infty,  \nonumber
\ee
\be h_2'(x_0)=0\Leftrightarrow x_0= \frac{1-\frac{\beta}{2}}{1-\frac{\alpha}{2}},  \nonumber
\ee
it is easy to see that $x_0$ is the maximum point of $h_2(x)$ in the interval $(0,\infty)$, so
\be\label{2cd2}  F_{22}(x)-F_{12}(x)<0\Leftrightarrow h_2(x)<0\Leftrightarrow \nu>(p-1)/d_2(\alpha, \beta).
\ee
Note that
 \be  F_{11}-F_{21}=x_1^{p-2}\{(p-1)-\nu\alpha(1-\frac{\alpha}{2})(\frac{x_2}{x_1})^{\frac{\beta}{2}}
-\frac{1}{2}\alpha\beta\nu(\frac{x_2}{x_1})^{\frac{\beta}{2}-1}\}.  \nonumber
\ee
Let
\be   h_3(x):= (p-1)-\nu\alpha(1-\frac{\alpha}{2})x^{\frac{\beta}{2}}
-\frac{1}{2}\alpha\beta\nu x^{\frac{\beta}{2}-1},\quad  0<x:=\frac{x_2}{x_1}<\infty,  \nonumber
\ee
then
\be h_3'(x_0)=0\Leftrightarrow x_0= \frac{1-\frac{\beta}{2}}{1-\frac{\alpha}{2}},  \nonumber
\ee
and it is easy to see that $x_0$ is the maximum point of $h_3(x)$ in the interval $(0,\infty)$. Therefore,
\be\label{2cd3}  F_{11}(x)-F_{21}(x)<0\Leftrightarrow h_2(x)<0\Leftrightarrow \nu>(p-1)/d_3(\alpha, \beta).
\ee
Combine with \eqref{2cd1}, \eqref{2cd2}, \eqref{2cd3}, we see that
\be  \sum_{j=1}^2F^{ij}>0,\;\;  i=1,2.  \nonumber
\ee
Then by Lemma \ref{2pringe}, the conclusion follows.
\ep

\bl\label{2deri} Consider the symmetric matrix $\gamma=(\gamma_{ij})$.   Assume  {\rm det}$(\gamma)\neq 0 $ and
\be\label{2temp}   \sum_{k=1}^r\gamma_{jk}c_k=1, \;\; j=1,...,r.
\ee
View  $(\sum_{j=1}^r c_j)$ as  a function of $\gamma_{ml}$. Then
\be  \frac{\partial}{\partial\gamma_{ml}}\Big(\sum_{j=1}^r c_j\Big)=-c_m c_l. \nonumber
\ee
\el

\bp
Let $(\gamma^{sj})$  represent the inverse matrix of $(\gamma_{jk})$.
Derivative on both side of \eqref{2temp} with respect to $\gamma_{ml}$ for any fixed $m,l$, we have
\be  \sum_k ( \gamma_{jk}c'_k+\delta_{mj}\delta_{kl}c_k) =0,  \nonumber \ee
that is,
\be\sum_k \gamma_{jk}c'_k=-\delta_{mj}c_l, \nonumber
\ee
\be \sum_j\sum_k \gamma^{sj} \gamma_{jk}c'_k=-\sum_j\delta_{mj}c_l\gamma^{sj},\nonumber
\ee
\be \sum_k\delta_{sk}c'_k=-\sum_j\delta_{mj}c_l\gamma^{sj},\nonumber
\ee
\be c'_s=c_l \gamma^{sm}.\nonumber
\ee
Hence,
\be\label{2temp2}  \sum_{s=1}^r c'_s=-c_l\sum_{s=1}^r \gamma^{sm}.
\ee
On the other hand, we see from \eqref{2temp} that
\be \sum_j\sum_k \gamma_{jk}c_k \gamma^{sj}=\sum_j \gamma^{sj}, \nonumber
\ee
\be \sum_k \delta_{ks}c_k=\sum_j \gamma^{sj}, \nonumber
\ee
hence,
\be c_s=\sum_j \gamma^{sj}. \nonumber
\ee
Since the matrix $\gamma=(\gamma_{ij})$ is symmetric, combine the equality above with \eqref{2temp2}, we have
\be    \sum_{s=1}^r c'_s=-c_m c_l. \nonumber
\ee
\ep

\noindent{\bf The proof of Theorem \ref{2th3}}.  Consider the matrix defined by
\be\label{zwm-aug-jzgm}\gamma:=
\begin{pmatrix}
       \gamma_{11} & \gamma_{12} &  \cdots & \gamma_{1r}  \\
        \gamma_{21} & \gamma_{22} &  \cdots & \gamma_{2r}  \\
        \vdots & \vdots &   \ddots &   \vdots \\
        \gamma_{r1} & \gamma_{r2} &  \cdots & \gamma_{rr}
\end{pmatrix},\ee
where $\gamma_{jk}=\gamma_{kj}$.
We consider the following critical elliptic system involving Hardy singular terms
\be\label{2model2}\begin{cases}-\Delta u_j-\frac{\lambda}{|x|^2} u_j
=\gamma_{jj}u_j^3+\sum\limits_{k\neq j}\gamma_{ij}u_i^2u_j,\quad x\in\R^4,\\
u_j(x)>0,\quad x\in \R^4\setminus\{0\}, \quad j=1,...,r.\end{cases}\ee

\br\label{zwm-aug-LM}  As previous definitions in  \eqref{2JJJ}, \eqref{zwm=110} and \eqref{zwm=111}, we may introduce the corresponding functional, Nehari manifold and the least energy  for the Eq. \eqref{2model2}. We adopt the same notations by $J, {\mathcal{N}}, \theta$ respectively as defined in \eqref{2JJJ}, \eqref{zwm=110} and \eqref{zwm=111} though the constants (coefficients) are replaced by those corresponding to   Eq. \eqref{2model2}. \er

Recall that  the matrix  $\gamma$ is invertible and the sum of each row  of the inverse    matrix  $\gamma^{-1}$ is greater than 0, it follows that the equation
\be\label{2w}    \sum_{k=1}^r \gamma_{jk} c_k=1,\;\;  j=1,...,r ,
\ee
has a solution $(c_1,...,c_r)$ satisfying  $c_j>0\; (j=1,...,r)$ and so $(\sqrt{c_1}z, ..., \sqrt{c_r}z)$ is a nontrivial solution of \eqref{2model2} (where $z$ is a solution of \eqref{zwm-aug-th13}) and
\be \label{2IJ-zwm}  \Theta=J(\sqrt{c_1}z, ..., \sqrt{c_r}z)=\sum_{j=1}^r c_j \Theta_1,\;  \hbox{ where } \Theta_1=I_\la(z) \;(see \; \eqref{zwm-aug-15} \; with\;
\lambda_j=\lambda).
\ee

\noindent{\bf The proof of Theorem\ref{2th3}-(1)}.   Let $\{(u_{1,n}, ..., u_{r,n})\}\subset \mathcal{N}$ be a minimizing sequence
 for $\Theta$, that is, $J(u_{1,n},...,u_{r,n})\rightarrow \Theta.$
Define
\be d_{i,n}=\Big(\int_{\Om}u_{i,n}^4 dx\Big)^{1/2}, i=1,...,r.  \nonumber\ee
Then by \eqref{2w}, we have
\begin{equation*}\aligned
2\sqrt{\Theta_1}d_{j,n}&\leq\int_{\R^4}|\nabla u_{j,n}|^2-\frac{\la}{|x|^2} u_{j,n}^2\\
&=\int_{\R^4}\gamma_{jj} u_{j,n}^4+\sum_{k\neq j}\int_{\R^4}\gamma_{kj}u_{k,n}^2 u_{j,n}^2\\
&\leq\gamma_{jj} d_{j,n}^2+\sum_{k\neq j}\gamma_{kj}d_{j,n}d_{k,n}.
\endaligned
\end{equation*}
On the other hand
\be 2\sqrt{\Theta_1}\sum_{i=1}^r d_{i,n}\leq 4 J(u_{1,n},...,u_{r,n})\leq 4\sum_{j=1}^r c_j \Theta_1+o(1),\nonumber \ee
thus we have
\be \begin{cases}\sum_{i=1}^r d_{i,n}\leq \sum_{i=1}^r c_i 2\sqrt{\Theta_1}+o(1),\\
\gamma_{ii}d_{i,n}+\sum_{k\neq j}\gamma_{ki}d_{k,n}\geq 2\sqrt{\Theta_1}.
  \end{cases}  \nonumber \ee
Recall   \eqref{2w}, then the inequalities above are equivalent to
\be \begin{cases}\sum_{i=1}^r (d_{i,n}- c_i 2\sqrt{\Theta_1})\leq o(1),\\
\gamma_{ii} (d_{i,n}- c_i 2\sqrt{\Theta_1})+\sum_{k\neq i}\gamma_{ki}(d_{k,n}-c_k 2\sqrt{\Theta_1})\geq0,\\
i=1,...,r.
\end{cases}
\nonumber\ee
By Lemma \ref{2principle}, we have $d_{i,n}\rightarrow c_i 2\sqrt{\Theta_1}$ as $n\rightarrow\infty$, and
\be 4\Theta=\lim_{n\rightarrow\infty}4J(u_{1,n} ,..., u_{r,n})\geq
 \lim_{n\rightarrow\infty}2\sqrt{\Theta_1}\sum_{i=1}^N d_{i, n}=4\sum_{i=1}^r c_i \Theta_1. \nonumber
 \ee
Combining this with \eqref{2IJ-zwm}, one has that
\be \Theta=\sum_{j=1}^r c_j \Theta_1=J(\sqrt{c_1}z, ..., \sqrt{c_r}z),  \nonumber \ee
and so $(\sqrt{c_1}z, ..., \sqrt{c_r}z)$ is a positive least energy solution of \eqref{2model2}.  \hfill$\Box$

\vskip0.3in

\noindent{\bf The proof of Theorem \ref{2th3}-(2)}. Namely,  we need to  prove the uniqueness of the ground state of \eqref{2model2}.
Let $(u_{1,0}, ..., u_{r,0})$ be any  least energy solution of \eqref{2model2}. Firstly we define the real functions with variables
$(t_1, ..., t_r)\in \R^r$:
\be\label{2fjfj}  f_j(t_1, ..., t_r):=\int_{\R^4}t_j\gamma_{jj} u_{j,0}^4+\sum_{k\neq j}\int_{\R^4}t_k\gamma_{kj}u_{k,0}^2 u_{j,0}^2
-\int_{\R^4}|\nabla u_{j,0}|^2-\frac{\la}{|x|^2} u_{j,0}^2.
\ee
Here we  regard   $\gamma_{ml}$ (for any fixed $(m,l)$ satisfying $  1\leq m, l\leq r$) as  the variable.
Recalling the definitions of$J, \mathcal{N}$ and $\Theta$, they all depend on $\gamma_{ml}$. Hence, we now adopt  the
 notations  $J({\gamma_{ml}})$,  $\mathcal{N}({\gamma_{ml}})$ and $\Theta({\gamma_{ml}})$  in this proof.
With the definitions above, we have $f_j(1,...,1)=0$  and
\be \frac{\pl f_j}{\pl t_i}=\gamma_{ij}\int_{\R^4}u_{i,0}^2 u_{j,0}^2.  \nonumber
\ee
Define the matrix:
\be F:= \Big(\frac{\pl f_j}{\pl t_i}\mid_{(1,...,1)}\Big). \nonumber \ee
Since the matrix $\gamma$ defined in \eqref{zwm-aug-jzgm} is positively definite, so is the following matrix $(\gamma_{ij}\int_{\R^4}u_{i,0}^2 u_{j,0}^2)$. Hence, $\det(F)>0$.
Therefore,  by the Implicit Function Theorem, the functions $t_j(\widetilde{\bb_{ml}})$ are well defined
and of class $C^1$ on $(\gamma_{ml}-\dd_1, \gamma_{ml}+\dd_1)$ for some $0<\dd_1\leq\dd$. Moreover, $t_j(\gamma_{ml})=1, j=1,...,r$,
and so we may assume that $t_j(\widetilde{\gamma_{ml}})>0$ for all $\widetilde{\gamma_{ml}}\in (\gamma_{ml}-\dd_1, \gamma_{ml}+\dd_1)$ by choosing a small $\dd_1$. From $f_k(t_1(\widetilde{\gamma_{ml}}), ..., t_r(\widetilde{\gamma_{ml}}))\equiv0$, it is easy to prove that:
\be  \sum_{j=1}^N\frac{\pl f_k}{\pl t_j}t_j'(\gamma_{ml})=-\frac{\pl f_k}{\pl \gamma_{ml}}.  \nonumber
\ee
Hence
\be t_j'(\gamma_{ml})=-\sum_{k=1}^N\frac{\pl f_k}{\pl \gamma_{ml}}\frac{F^*_{kj}}{\det(F)}, \nonumber
\ee
here $F^*:=(F^*_{kj})$ denotes the adjoint matrix of $F$. From \eqref{2fjfj}, we have
\be \frac{\pl f_k}{\pl \gamma_{ml}}=\dd_{km}\int_{\R^4}u_{m,0}^2 u_{l,0}^2 dx, \;\; \hbox{ where } \delta_{km} \hbox{ is the Kronecker notation},\nonumber\ee
hence
\be t_j'(\gamma_{ml})=-\sum_k \dd_{km}\int_{\R^4}u_{m,0}^2 u_{l,0}^2 dx \frac{F^*_{kj}}{\det(F)},  \nonumber \ee
that is,
\be\label{2tm}
t_j'(\gamma_{ml})=-\int_{\R^4}u_{m,0}^2 u_{l,0}^2 dx \frac{F^*_{mj}}{\det(F)}.
\ee
By  the Taylor's expansion, we see that
\be t_j'(\widetilde{\gamma_{ml}})=1+t_j'(\gamma_{ml})(\widetilde{\gamma_{ml}}-\gamma_{ml})+O((\widetilde{\gamma_{ml}}-\gamma_{ml})^2).\nonumber \ee
Note that $t_j(t_1(\widetilde{\gamma_{ml}}) ,..., t_r(\widetilde{\gamma_{ml}}))\equiv0$ implies that
$$ (\sqrt{t_1({\widetilde{\gamma_{ml}}})}u_{1,0}, ..., \sqrt{t_1({\widetilde{\gamma_{ml}}})}u_{r,0})\in \mathcal{N}({\widetilde{\gamma_{ml}}}),$$
therefore
\be \aligned J(\widetilde{\gamma_{ml}})&\leq E_{\widetilde{\gamma_{ml}}}(\sqrt{t_1({\widetilde{\gamma_{ml}}})}u_{1,0},...,\sqrt{t_1({\widetilde{\gamma_{ml}}})}u_{r,0})\\
&=\frac{1}{4}\sum_{j=1}^r t_j(\widetilde{\gamma_{ml}}\int_{\R^4}|\nabla u_{j,0}|^2-\frac{\la}{|x|^2} u_{j,0}^2dx\\
&=J(\gamma_{ml})+\frac{1}{4}D(\widetilde{\gamma_{ml}}-\gamma_{ml})+O((\widetilde{\gamma_{ml}}-\gamma_{ml})^2), \nonumber
\endaligned
\ee
where
\be \aligned
D:&=\sum_{j=1}^r t_j'(\widetilde{\gamma_{ml}})\int_{\R^4}|\nabla u_{j,0}|^2-\frac{\la}{|x|^2} u_{j,0}^2dx\\
&=\sum_{j=1}^r t_j'(\widetilde{\gamma_{ml}})\Big(\int_{\R^4}\gamma_{jj} u_{j,0}^4dx+\sum_{k\neq j}\gamma_{kj}\int_{\R^4}u_{k,0}^2 u_{j,0}^2 dx\Big)\\
&=-\int_{\R^4}u_{m,0}^2 u_{l,0}^2 \sum_{j=1}^r\frac{F_{mj}^*}{det(F)}(F_{jj}+\sum_{k\neq j}F_{kj})\\
&=-\int_{\R^4}u_{m,0}^2 u_{l,0}^2\frac{1}{det(F)}\sum_{k=1}\sum_{j=1}F_{mj}^*F_{kj}\\
&=-\int_{\R^4}u_{m,0}^2 u_{l,0}^2\frac{1}{det(F)}\sum_{k=1}\dd_{km}det(F)\\
&=-\int_{\R^4}u_{m,0}^2 u_{l,0}^2. \nonumber
\endaligned
\ee
Here we have used \eqref{2tm}.
It follows that
\be\frac{ J(\widetilde{\gamma_{ml}})- J(\gamma_{ml})}{\widetilde{\gamma_{ml}}-\gamma_{ml}}\geq \frac{D}{4}+O(\widetilde{\gamma_{ml}}-\gamma_{ml})\nonumber
\ee
as $\widetilde{\gamma_{ml}}\nearrow \gamma_{ml}$ and so $J'(\gamma_{ml})\geq \frac{D}{4}$.
Similarly, we have $\frac{ J(\widetilde{\gamma_{ml}})- J(\gamma_{ml})}{\widetilde{\gamma_{ml}}-\gamma_{ml}}\leq\frac{D}{4}
+O(\widetilde{\gamma_{ml}}-\gamma_{ml})$ as $\widetilde{\gamma_{ml}}\searrow \gamma_{ml}$, that is,  $J'(\gamma_{ml})\leq\frac{D}{4} $ in this case.
Hence,
\be J'(\gamma_{ml})= \frac{D}{4}=-\frac{1}{4}\int_{\R^4}u_{m,0}^2 u_{l,0}^2. \nonumber
\ee
On the other hand, by  Lemma \ref{2deri}, we have
\be  J'(\gamma_{ml})=-c_m c_l \Theta_1=-\frac{1}{4}c_m c_l\int_{\R^4}z^4. \nonumber \ee
Hence,
\be \int_{\R^4}u_{m,0}^2 u_{l,0}^2=c_m c_l\int_{\R^4}z^4. \nonumber
\ee
Define
\be(\widetilde{u}_{1}, ..., \widetilde{u}_{r}):=(\frac{1}{\sqrt{c_1}}u_{1,0} ,..., \frac{1}{\sqrt{c_r}}u_{r,0}). \nonumber \ee
Combine this with the following identity:
\begin{equation*}
\gamma_{jj}.c_j+\Sigma_{k\neq j}\beta_{kj}.c_k=1
\end{equation*} and $(u_{1,0},...,u_{r,0})\in\mathcal{N}$, we get that
\be
\aligned
\int_{\R^4}(|\nabla \widetilde{u}_j|^2-\frac{\la}{|x|^2} \widetilde{u}_j^2)dx
&=\frac{1}{{c_j}}\int_{\R^4}(|\nabla u_{j,0}|^2-\frac{\la}{|x|^2} u_{j,0}^2)dx\\
&=\frac{1}{{c_j}}\Big(\int_{\R^4}\gamma_{jj} u_{j,0}^4dx+\sum_{k\neq j}\gamma_{kj}\int_{\R^4}u_{k,0}^2 u_{j,0}^2 dx\Big)\\
&=\frac{1}{{c_j}}\Big(\mu_j c_j^2+\sum_{k\neq j}\gamma_{kj}c_k c_j\Big)\int_{\R^4}\omega^4dx\\
&=\Big(\gamma_{jj} c_j+\sum_{k\neq j}\gamma_{kj}c_k\Big)\int_{\R^4}\omega^4dx\\
&=\int_{\R^4}\widetilde{u}_j^4dx. \nonumber
\endaligned
\ee
Then by \eqref{2w}, we have
\begin{equation}
\frac{1}{4}\int_{\R^4}\Big(|\nabla \widetilde{u}_j|^2-\frac{\la}{|x|^2} \widetilde{u}_j^2 \Big)dx\geq \Theta_1,\;\;\; j=1,...,r.
\end{equation}
Hence,
\begin{equation*}
\aligned
\Theta &=\sum_{j=1}^r c_j \Theta_1=\frac{1}{4}\sum_{j=1}^r\int_{\R^4}(|\nabla {u}_{j,0}|^2-\frac{\la}{|x|^2} {u}_{j,0}^2 )dx\\
&=\frac{1}{4}\sum_{j=1}^r c_j\int_{\R^4}(|\nabla \widetilde{u}_j|^2-\frac{\la}{|x|^2} \widetilde{u}_j^2)dx\\
&\geq \sum_{j=1}^r c_j \Theta_1.
\endaligned
\end{equation*}
This implies that
\begin{equation*}
\aligned
\frac{1}{4}\int_{\R^4}(|\nabla \widetilde{u}_{j}|^2-\frac{\la }{|x|^2}\widetilde{u}_{j}^2 )dx=\Theta_1,\;\;j=1,...,r.
\endaligned
\end{equation*}
Then we see that $\widetilde{u}_j (j=1,...,r)$ are the positive least energy solutions of \eqref{2hardyeq}.
We see from the fact that
\begin{equation*}
-\Delta \widetilde{u}_j-\frac{\la}{|x|^2} \widetilde{u}_j=\gamma_{jj} c_j \widetilde{u}_j^3+\sum_{k\neq j}\gamma_{kj}c_k \widetilde{u}_k^2 \widetilde{u}_j=\widetilde{u}_j^3
\end{equation*}
and
\begin{equation*}
\gamma_{jj} c_j \widetilde{u}_j^2+\sum_{k\neq j}\gamma_{kj}c_k \widetilde{u}_k^2 =\widetilde{u}_j^2,
\end{equation*}
hence
\begin{equation*}
\gamma_{jj} c_j +\sum_{k\neq j}\gamma_{kj}c_k (\frac{\widetilde{u}_k}{\widetilde{u}_j})^2 =1.
\end{equation*}
Since the matrix $\gamma$ is invertible, we get that
$\frac{\widetilde{u}_k}{\widetilde{u}_j}=1, k\neq j$.  That is, $\widetilde{u}_k=\widetilde{u}_j, k\neq j$.
Denote that $U=\widetilde{u}_1$, then $(u_{1,0},...,u_{r,0})=(\sqrt{c_1}U,..., \sqrt{c_r}U)$,
where $U$ is a positive least energy solution of \eqref{2hardyeq}.\hfill$\Box$

\vskip0.3in

\noindent{\bf The proof of Theorem \ref{2th4}}.
We consider the following doubly critical Shr\"{o}dinger system (i.e., \eqref{2mod4}) on $\R^N$:
\begin{equation}\label{2th4mo}
\begin{cases}
-\Delta u-\frac{\lambda}{|x|^2}u=u^{2^*-1}+\nu \alpha u^{\alpha-1}v^{\beta}, \\
-\Delta v-\frac{\lambda}{|x|^2}v=v^{2^*-1}+\nu \alpha u^{\alpha}v^{\beta-1}.
\end{cases}
\end{equation}
Let $p=2^\ast/2.$ It is easy to see that the following system
\begin{equation}\label{2f1f2}
\begin{cases}
f_1(x_1,x_2):&=x_1^{p-1}+\nu\alpha x_1^{\frac{\alpha}{2}-1}x_2^{\frac{\beta}{2}}=1, \\
f_2(x_1,x_2):&=x_2^{p-1}+\nu\beta x_1^{\frac{\alpha}{2}}x_2^{\frac{\beta}{2}-1}=1,
\end{cases}
\end{equation}
admits a positive solution $(c_1, c_2)$ for any $\nu>0$. In fact, from the first equality of \eqref{2f1f2}, we know  that
$x_2=(2\nu\alpha)^{-\frac{2}{\beta}}(1-x_1^{p-1})^{\frac{2}{\beta}}x_1^{\frac{2-\alpha}{\beta}}$. The system \eqref{2f1f2} admitting
 a positive solution is equivalent to the equation
\be\aligned f(x_1):=&(\nu\alpha)^{-\frac{1}{\beta}(p-1)}x_1^{\frac{2-\alpha}{\beta}(p-1)}(1-x_1^{p-1})^{\frac{2}{\beta}(p-1)}\\
&+\nu\beta(\nu\alpha)^{-\frac{\beta-2}{\beta}}x_1^{\frac{\alpha+\beta-2}{\beta}}(1-x_1^{p-1})^{-\frac{2-\beta}{\beta}}-1=0  \nonumber
\endaligned\ee
has a root  in the interval $(0,1)$. Since $f(0)=-1$ and $\lim_{x_1\rightarrow 1^{-}}f(x_1)=+\infty$, the conclusion
follows from the Mean Value Theorem.
Hence, $(\sqrt{c_1}z, \sqrt{c_2}z)$ is a nontrivial solution of \eqref{2th4mo} and
\be\label{2t4tem}  0<\Theta\leq J(\sqrt{c_1}z, \sqrt{c_2}z)=(c_1+c_2)\Theta_1.
\ee
Now we assume that $\nu>(p-1)/\min\{d_1(\alpha, \beta),d_2(\alpha, \beta),d_3(\alpha, \beta)\}$, and we shall prove that
$\Theta=J(\sqrt{c_1}z,\sqrt{c_2}z)$. Let $\{(u_n,v_n)\}\subset \mathcal{N}$ be a minimizing sequence for $\Theta$,
that is,  $J(u_n,v_n)\rightarrow \Theta$. Define
\be d_{1,n}=\Big(\int_{\R^N}|u_n|^{2p}dx\Big)^{\frac{1}{p}},\;\;  d_{2,n}=\Big(\int_{\R^N}|v_n|^{2p}dx\Big)^{\frac{1}{p}}.\nonumber
\ee
By \eqref{zwm-aug-15} and \eqref{2JJJ}, we have
 \begin{equation*}
 \aligned
 (N\Theta_1)^{2/N}d_{1,n}\leq \int_{\R^N}(|\nabla u_n|^2-\frac{\la}{|x|^2}u_n^2)&=\int_{\R^N}(|u_n|^{2p}+\nu\alpha|u_n|^{\alpha}v_n^{\beta}) \\
  &\leq d_{1,n}^p+\nu\alpha d_{1,n}^{\alpha/2}d_{2,n}^{\beta/2},\\
  (N\Theta_1)^{2/N}d_{2,n}\leq \int_{\R^N}(|\nabla v_n|^2-\frac{\la}{|x|^2}v_n^2)&=\int_{\R^N}(|v_n|^{2p}+\nu\beta|u_n|^{\alpha}v_n^{\beta}) \\
 &\leq d_{2,n}^p+\nu\beta d_{1,n}^{\alpha/2}d_{2,n}^{\beta/2}.
 \endaligned\end{equation*}
Since $J(u_n,v_n)=\frac{1}{N}\int_{\R^N}(|\nabla u_n|^2+|\nabla v_n|^2-\frac{\la}{|x|^2}u_n^2-\frac{\la}{|x|^2}v_n^2)$,
by \eqref{2t4tem}, we have
 \begin{equation}\label{2t4p}\begin{cases}
 \aligned
  (N\Theta_1)^{2/N}(d_{1,n}+d_{2,n})\leq N J(u_n,v_n)&\leq N(c_1+c_2)\Theta_1+o(1),\\
 d_{1,n}^{p-1}+\nu\alpha d_{1,n}^{\alpha/2-1}d_{2,n}^{\beta/2}&\geq (N\Theta_1)^{2/N},\\
 d_{2,n}^{p-1}+\nu\beta d_{1,n}^{\alpha/2}d_{2,n}^{\beta/2-1}&\geq (N\Theta_1)^{2/N}.
 \endaligned\end{cases}
 \end{equation}
First, this means that  $d_{1,n},d_{2,n}$  are uniformly bounded.
Passing to a subsequence we may assume that $d_{1,n}\rightarrow d_1,d_{2,n}\rightarrow d_2$. It is easy to check that  $d_1>0, d_2>0$.
Denote
\be x_1=\frac{d_1}{(N\Theta_1)^{1-\frac{N}{2}}},\;\;  x_2=\frac{d_2}{(N\Theta_1)^{1-\frac{N}{2}}}. \nonumber\ee
By a simple scaling we can transform  \eqref{2t4p}   to
 \begin{equation*}\begin{cases}
 \aligned
x_1+x_2\leq c_1+c_2,\\
 x_1^{p-1}+\nu\alpha x_1^{\alpha/2-1}x_2^{\beta/2}&\geq 1,\\
 x_2^{p-1}+\nu\beta x_1^{\alpha/2}x_2^{\beta/2-1}&\geq 1.
\endaligned\end{cases}\end{equation*}
By Lemma \ref{2corrth4}, we see that   $x_1=c_1, x_2=c_2$. It follows that
\be d_{1,n}\rightarrow c_1(N\Theta_1)^{1-\frac{N}{2}}, d_{2,n}\rightarrow c_2(N\Theta_1)^{1-\frac{N}{2}},\;\;\hbox{as}\;\; n\rightarrow \infty\nonumber
\ee
and
\be  N\Theta=\lim_{n\rightarrow\infty}NJ(u_n, v_n)\geq(N\Theta_1)^{2/N}(d_{1,n}+d_{2,n})=N(c_1+c_2) \Theta_1. \nonumber
\ee
Combing this with  \eqref{2t4tem}, we have
\be\Theta=(c_1+c_2)\Theta_1=J(\sqrt{c_1}z, \sqrt{c_2}z), \nonumber
\ee
and therefore, $(\sqrt{c_1}z, \sqrt{c_2}z)$ is a positive least energy solution of \eqref{2th4mo}. \hfill $\Box$

\vskip0.36in





\end{document}